\documentclass{amsart}

\usepackage{amssymb}
\usepackage{amsfonts}
\usepackage{epsfig}
\input{arrow}

\copyrightinfo{}{}
\issueinfo{}
  {}
  {}
  {}

\PII{}
\def\@serieslogo{\@empty}

\newcommand{\bbold}{\mathbb}

\newcommand{\cal}{\mathcal}

\newcommand{\rom}{\textup}

\def\tm{\subseteq}
\def\om{\supseteq}

\def\id{\operatorname{id}}
\def\ohne{\setminus}
\def\leer{\emptyset}

\def\root{\operatorname{root}}
\def\paths{\operatorname{paths}}

\renewcommand{\leq}{\leqslant}
\renewcommand{\geq}{\geqslant}
\renewcommand{\epsilon}{\varepsilon}

\def\iff{\Longleftrightarrow}

\def\ini{\preceq}
\def\propini{\prec}
\def\tri{\vartriangleleft}

\def \R { {\bbold R} }
\def \Q { {\bbold Q} }
\def \Z { {\bbold Z} }
\def \N { {\bbold N} }

\newtheorem{theorem}{Theorem}[section]
\newtheorem*{theorem-unnumbered}{Theorem}
\newtheorem{lemma}[theorem]{Lemma}
\newtheorem{prop}[theorem]{Proposition}
\newtheorem{cor}[theorem]{Corollary}

\theoremstyle{definition}
\newtheorem{definition}[theorem]{Definition}
\newtheorem{algorithm}[theorem]{Algorithm}

\theoremstyle{remark}
\newtheorem*{example}{Example}
\newtheorem*{examples}{Examples}
\newtheorem{examplesNumbered}[theorem]{Examples}
\newtheorem*{notations}{Notations}
\newtheorem*{remarks}{Remarks}
\newtheorem*{remark}{Remark}

\def \E { {\bbold E} }
\newcommand{\boproof}{\begin{proof}}
\newcommand{\eoproof}{\end{proof}}
\newcommand{\red}{\sqsubseteq}

\newcommand{\Kw}[1]{\underline{\bf #1}}
\DeclareMathOperator{\normalForm}{normalForm}

\numberwithin{equation}{section}

\newcommand{\abs}[1]{\lvert#1\rvert}

\begin{document}
\title{Finiteness Theorems in Stochastic Integer Programming}

\dedicatory{Dedicated to the memory of C.~St.~J.~A. Nash-Williams, 1932--2001.}

\author{Matthias Aschenbrenner}
\address{Department of Mathematics, Statistics, and Computer Science\\
University of Illinois at Chicago\\
815 S. Morgan St. (M/C 249) \\
Chicago, IL 60607-7045\\
U.S.A.}
\email{maschenb@math.uic.edu}

\author{Raymond Hemmecke}
\address{Institut f\"ur Mathematische Optimierung\\
Fakult\"at f\"ur Mathematik \\
Universit\"at Magdeburg\\
Universit\"atsplatz 2\\
39106 Magdeburg \\
Germany}
\email{raymond@hemmecke.de}
\thanks{The work of the first author was partially supported by NSF grant
DMS 03-03618. The work on this paper was begun in 2001, while the first author
was a Charles B.~Morrey Jr.~Visiting Research Professor at UC Berkeley and the second author was a
Visiting Research Assistant Professor at UC Davis. The support of these institutions is gratefully
acknowledged. The
authors would also like to thank Bernd Sturmfels for his
encouragement to write this paper.}
\date{January 2005.}



\begin{abstract}
We study Graver test sets for families of linear multi-stage stochastic integer
programs with varying number of scenarios.
We show that these test sets can be decomposed into finitely many
``building blocks'', independent of the number of scenarios,
and we give an effective procedure to compute these
building blocks. The paper includes an introduction to
Nash-Williams' theory of better-quasi-orderings, which is used to show
termination of our algorithm. We also apply
this theory to finiteness results for Hilbert functions.
\end{abstract}

\maketitle

\setcounter{tocdepth}{1}
\tableofcontents

\newpage

\section*{Introduction}

\noindent
A wide range of practical optimization problems can be modeled as
(mixed-integer) linear programs. While some methods for solving problems
like these assume
all data to be known in advance, others deal with uncertainty. Such
uncertainty is often inherent to the problem, e.~g., demands, prices,
or network structure. To be able to treat optimization
problems of this type, one requires  stochastic information about the uncertain
components. With this knowledge, the problem can then
be formulated as a {\em stochastic}\/ linear program: for example,
one that minimizes {\em expected}\/ costs, or one that minimizes
the {\em risk}\/ that the costs exceed a given threshold.

A typical situation is the following: In the first step, one has to make a decision $x$ without knowing the outcome
of a random event that lies in the future. For example, one needs to
decide where to open new production facilities without exactly knowing  the
demand in each geographical area. In a second stage, after
observing the uncertain data (here: demands), one can make a second
(recourse/repair) decision $y$. For example, one may relocate
facilities, decide on production or transportation plans, etc.
The goal
is now to maximize revenues or profits from these facilities. These
quantities can be calculated from the immediate costs $c^\intercal x$
for opening a new location (first stage decision) plus the expected
costs and revenues ${\cal Q}(x)$ from relocating and running the
facilities (second stage decision).
This  is a typical instance of an {\em expectation problem,}\/ which can
be modeled as
\[
\min\big\{c^\intercal x+{\cal Q}(x):Ax=b,\ x\in\N^m\big\},
\]
where  $A\in\Z^{l\times m}$, $b\in\Z^l$, $c\in\R^m$, and
\[
{\cal Q}(x)=\E_\omega\min\big\{d(\omega)^\intercal y: Wy=h(\omega)-Tx,
\ y\in\N^n\big\}
\]
where $T\in\Z^{l\times m}$, $W\in\Z^{l\times n}$, and $d$, $h$ 
are  random
variables taking values in $\R^n$ and $\Z^l$, respectively.
Note that in the  formulation of the problem we assume that the outcome of
the random event $\omega$ does not depend on the first-stage decision $x$
that has been made.

When no integrality constraints are imposed on $y$, the
computation of ${\cal Q}(x)$ involves the symbolic computation of
an integral, which is often very hard (even for relatively few variables).
In practice, the computation can sometimes be simplified by exploiting
continuity and convexity of ${\cal Q}(x)$.
However, if $y$ is required to be integer, ${\cal
Q}(x)$ is only semi-continuous, and usually not convex.
Hence, to make computations possible, the probability distribution for
$\omega$ is usually approximated by $N$ {\em scenarios}\/ $\omega_1,\dots,\omega_N$
with respective probabilities
$\pi_1,\dots,\pi_N$. In this way,
the integral involved in the computation of ${\cal Q}(x)$ becomes
a sum, and we obtain an integer program having
a separate set of $y$-variables corresponding to each scenario:
\[
\min\left\{c^\intercal x+\sum_{i=1}^N \pi_i\cdot(d^\intercal y_i):Ax=a,\
Wy_i=h_i-Tx,\ x\in\N^m,\ y_i\in\N^n\right\}.
\]
The coefficient matrix of this problem
has a nice block structure:
\[
\left(
\begin{array}{cccc}
A & & & \\
T & W & & \\
\vdots & & \ddots & \\
T &   & & W \\
\end{array}
\right)
\]
This two-stage setup (making a decision, observing the outcome of
the random event, and then making a recourse decision) can be
iterated, leading to the notion of a {\em multi-stage stochastic
integer program.}\/ (See \cite{Schultz+Roemisch} for a recent survey
of this topic.) In these programs,  information is revealed only at
certain points in time, and decisions have to be made without
knowing the outcome of future random events. Again we assume that
each random event neither depends on the outcome of previous random
events nor on the decisions made in the previous stages.  The
integer programs which we obtain by discretizing the probability
distributions of the random variables involved quickly become very
big and hard to solve in practice. However, as in the two-stage
example above, the non-anticipativity assumption leads to highly
structured problem matrices.

In \cite{Hemmecke+Schultz:dec2SIP}, Hemmecke and Schultz exploited this
structure to construct a novel algorithm for the solution of
two-stage stochastic integer programs, which is based on successive
augmentation of a given feasible solution.
The main goal of this paper is to extend the ideas presented in
\cite{Hemmecke+Schultz:dec2SIP,Hemmecke:Diss} from
two-stage to multi-stage programs.
Instead of decomposing the problem itself, as essentially all
other decomposition approaches to two- and multi-stage stochastic integer
programming do, we will decompose an object that is closely related to
the given problem matrix: its
{\em Graver test set}\/ (or {\em Graver basis}\/).
The structure of the
problem matrix also imposes a lot of structure on the
elements in the Graver basis: for a given family
of $(k+1)$-stage stochastic integer programs ($k\in\N$)
we can define
a certain ({\it a priori}\/ infinite)
set ${\cal H}_{k,\infty}$ of ``building blocks,'' from which all
elements in the Graver basis of the coefficient matrix can be reconstructed,
{\em independent}\/ of the number $N$ of scenarios.
(We refer to Section~\ref{Building Blocks} for the precise definition.)
We will prove that this set
of building blocks, in fact, is always finite. We also
show how to compute ${\cal H}_{k,\infty}$
(even if only theoretically), and how it can be employed
to solve any given particular instance of the given family of
$(k+1)$-stage stochastic integer programs.


It is perhaps remarkable that
our proof of the finiteness of
$\cal H_{k,\infty}$ rests on some non-trivial properties of
the set $\N^n$ of $n$-tuples of natural numbers.
Some of those have appeared before,  under various guises,
in computational  algebra. The most prominent one, known as
``Dickson's Lemma,'' can  formulated like this:

\begin{quote}
{\it For every infinite sequence $X^{\nu^{(1)}},X^{\nu^{(2)}},\dots$
of monomials in the polynomial ring $\Q[X]=\Q[X_1,\dots,X_n]$,
where $X^\nu=X_1^{\nu_1}\cdots X_n^{\nu_n}$ for
$\nu=(\nu_1,\dots,\nu_n)\in\N^n$, there exist indices $i<j$ such that
$X^{\nu^{(i)}}$ divides $X^{\nu^{(j)}}$.}\/
\end{quote}

\noindent
This simple combinatorial
fact alone is at the heart of many finiteness phenomena
in commutative algebra,
since it readily implies Hilbert's Basis Theorem:
every ideal of the polynomial ring $K[X]=K[X_1,\dots,X_n]$, where $K$ is a
field, is finitely generated.
Dickson's Lemma yields that every monomial ideal of $K[X]$
(that is, an ideal generated by monomials) is finitely generated
(see, e.g.,  Proposition~\ref{Folk} below), and
Gordan's famous proof \cite{Gordan}
of the Hilbert Basis Theorem
extends this to all ideals of $K[X]$.
Now Dickson's Lemma in turn is a consequence of a somewhat
more powerful (and less obvious) finiteness statement:

\begin{quote}
{\it There is no infinite sequence $I^{(1)},I^{(2)},\dots$
of monomial ideals in the polynomial ring $\Q[X]$
such that $I^{(i)}\not\supseteq I^{(j)}$ for all $i<j$.}\/
\end{quote}

\noindent
This principle can be
easily shown using known techniques in the subject of
``well-quasi-orderings''.
Maclagan \cite{Maclagan} rediscovered this fact
(by primary decomposition of monomial ideals) and
demonstrated how it can be used to give short proofs of
several other finiteness statements like the
existence of universal Gr\"ob\-ner bases
and the finiteness of the number
of atomic fibers of a matrix with non-negative integer entries.
(Other applications can be found in \cite{Haiman-Sturmfels},
\cite{Hemmecke-atomic},
\cite{MaclaganThomas}.)
The connections to well-quasi-orderings have been made explicit and
further explored by the first author and Pong in
\cite{AschenbrennerPong}.

In this paper we exploit an  infinite hierarchy of finiteness principles, of
which the statements above only represent the two bottom  levels.
For a precise formulation we refer to
Theorem~\ref{NmBetterNoeth} below. Here, let us just state
an attractive consequence of the next statement in the hierarchy:

\begin{theorem-unnumbered}
Let $\cal S$ be a collection of monomial ideals in
the polynomial ring $\Q[X]$, and let
${\cal M}_1,{\cal M}_2,\dots$ be an infinite sequence
of collections of monomial
ideals from ${\cal S}$, where each $\cal M_i$ is closed under inclusion
\textup{(}if $I\in {\cal M}_i$ and $J\in\cal S$ is a monomial ideal such that
$J\subseteq I$, then $J\in {\cal M}_i$\textup{)}. Then
$\cal M_i\subseteq\cal M_j$ for some indices $i\neq j$.
\end{theorem-unnumbered}

(The statement remains true if ``closed under inclusion''
is replaced by ``closed under reverse inclusion.'')

We construct this hierarchy  using Nash-Williams' beautiful
theory of ``better-quasi-ordered sets''. This theory, although of
a fundamentally combinatorial nature, is probably less well-known
in the field of algorithmic algebra than among
logicians, who explored its connections to descriptive set theory
and computability theory  \cite{Simpson}, and, more
recently, investigated its logical strength \cite{Cholak-Marcone-Solomon},
\cite{Marcone}. It is for this
reason that we include an introduction to this subject in Part~1 of the
paper (Sections~\ref{Section-Preliminaries}--\ref{Section-Better-Noetherian}),
with the hope that it will become useful
as a general guide for proving termination of algebraic algorithms.
We finish Part~1 by applying Nash-Williams' theory to prove a few
finiteness properties for
Hilbert functions (in Section~\ref{Applications to Hilbert Functions}),
some of which are known (Corollary~\ref{HilbertCor}), and some of which might
be new (Proposition~\ref{H-is-Noetherian}).

In Part~2 of the paper we then apply the theorems of Part~1 to
establish finiteness and computability in the decomposition approach to
solve multi-stage stochastic integer programs mentioned above.
We begin by giving a brief introduction to Graver test sets
for integer linear programs. For a more thorough treatment see, e.g.,
\cite{Hemmecke:PSP}.
After introducing ${\cal H}_{k,\infty}$ in Section~\ref{Building Blocks},
we show that this set is
always finite, and give an algorithm to compute a set of vectors containing
it (in
Section~\ref{Computation of Building Blocks}).
The set ${\cal H}_{k,\infty}$ holds an enormous
amount of information. We finish the paper by
showing how knowledge of (an object related to) $\cal H_{k,\infty}$ allows
one to solve any given instance of our
family of $(k+1)$-stage stochastic integer programs, for any given
number of scenarios.

\subsection*{Conventions and notations}
Throughout this paper, $m$ and $n$ range over the set
$\N=\{0,1,2,\dots\}$ of
natural numbers.

\part{Noetherian Orderings and Monomial Ideals}

\section{Preliminaries}\label{Section-Preliminaries}

\noindent
We first introduce some notations 
about sets of natural numbers which are constantly used 
throughout Part~1. We then discuss an infinitary version of
Ramsey's Theorem due to Galvin and Prikry,
which is at the base of Nash-Williams' theory, together with the
notion of a ``barrier''. The reader may skip this section at first reading
and come back to it when it is really
needed (in Section~\ref{Section-Better-Noetherian}).

\subsection*{Sets and sequences of natural numbers}

For a set $X$, we denote the power set of $X$ (the set
of all subsets of $X$)  by ${\cal P}(X)$.
Given $X\subseteq\N$
and $n$, we denote by $[X]^n$ the set of
all subsets of $X$ consisting of exactly $n$ elements. We let
$$[X]^{<\omega}:=\bigcup_{n\in\N}\ [X]^n,$$ the set of all finite subsets
of $X$. By $[X]^\omega$ we denote the set of all infinite subsets
of $X$.
(So ${\cal P}(X)=[X]^{<\omega}\cup [X]^{\omega}$.)

For every subset of $\N$, there is a unique sequence enumerating it in
strictly increasing order; throughout this section, we will identify subsets of $\N$
and strictly increasing sequences of natural numbers in this way.
(The empty set $\leer$ corresponds to the empty sequence.)
Below, $s,t,u$ range over $[\N]^{<\omega}$, and $U,V,W,X$ over
$[\N]^\omega$.

We denote by 
$l(s)$ the cardinality of $s$.
So if $s$ is identified with the corresponding strictly increasing sequence,
then $l(s)$ is its {\bf length.} 
For every $0\leq i<l(s)$ we write $s_i$ for the $(i+1)$-st
element of $s$; therefore, we can write $s$ as $s=(s_0,\dots,s_{l(s)-1})$,
with $s_0<s_1<\cdots<s_{l(s)-1}$. 
We write $s \ini t$ if $s$ is an
initial segment of $t$, that is, $l(s)\leq l(t)$ and $s_i=t_i$ for all
$0\leq i<l(s)$. We put $s \propini t$ if $s$ is a proper initial segment
of $t$, i.e., $s\ini t$ and $s\neq t$. 
These relations extend in a natural way also to the case where $t$ 
is replaced by an infinite subset of $\N$.
Clearly $s\ini t$ implies $s\tm t$.

If $l(s)=n\geq 1$, then $s\ohne\{\min s\}$ is the
sequence $(s_1,\dots,s_{n-1})$ obtained from $s=(s_0,s_1,\dots,s_{n-1})$ 
by leaving out its first element $s_0$. 
Similarly if $l(s)\geq 1$, then
$t\ohne\{\max t\}$ is the sequence obtained from $t$ by
leaving out its last element.
We define, for non-empty $s,t$:
$$
s\tri t \qquad :\iff \qquad s\ohne\{\min s\}\ini t\ohne\{\max t\}.
$$
Note that the relation $\propini$ 
on $[\N]^{<\omega}$ is irreflexive ($s\not\propini s$ for all $s\in
[\N]^{<\omega}$) and  transitive (if 
$s\propini t$ and $t\propini u$, then $s\propini u$), whereas
the relation $\tri$ on the set $[\N]^{<\omega}$ 
is neither irreflexive nor transitive.

\subsection*{The Ellentuck topology}
We write
$s<U$ if $\max s < \min U$.
Here and below, $\max\emptyset:=-\infty<a$ for all $a\in\N$.
If $s<U$, we put
$$[s,U]:=\bigl\{X\in [s\cup U]^\omega : s\prec X\bigr\}.$$
We endow $[\N]^\omega$ with
the {\bf Ellentuck topology,} whose basic open sets are the sets of the form
$[s,U]$ for $s<U$ as above. 
(This topology was first introduced in \cite{Ellentuck}.)
For $X\in [\N]^\omega$, we consider each $[X]^\omega$ as a 
subspace of $[\N]^\omega$, equipped with the induced topology.


\subsection*{A theorem of Galvin and Prikry}

If $\N=P\cup Q$ is a partition of $\N$, then one of $P$ or $Q$ is infinite,
by the familiar (Dirichlet) pigeon-hole principle. Ramsey's Theorem
\cite{Ramsey}
is an extension of this principle: if $[\N]^n=P\cup Q$ is a
partition of $[\N]^n$, then there is $H\in [\N]^\omega$ such that
$[H]^n\tm P$ or $[H]^n\tm Q$. (Such a set $H$ is called a {\bf homogeneous}
set for the partition $[\N]^n=P\cup Q$.)
Later on, we will need a far-reaching generalization of this theorem, 
concerning partitions $[\N]^\omega=P\cup Q$ of $[\N]^\omega$. We have to
place some restrictions on the nature of the partitioning sets $P$ and $Q$,
since the natural analogue of Ramsey's Theorem for partitions of $[\N]^\omega$
fails for pathological partitions constructed using the Axiom of
Choice. (See \cite{Kechris}, Section~19.) 

\begin{theorem}\label{GPE} \rom{(Galvin-Prikry \cite{GalvinPrikry}.)}
Let $X\in [\N]^\omega$, and suppose $[X]^\omega=P\cup Q$ 
is a partition
of $[X]^\omega$, where $P$ is an open set 
\rom{(}in the Ellentuck topology\rom{)}. Then there
exists $H\in [X]^\omega$ such that $[H]^\omega\tm P$ or
$[H]^\omega\tm Q$.
\end{theorem}

We refer to \cite{Kechris}, Section~19 or \cite{Simpson} for
a proof. 

\subsection*{Blocks and barriers}

The language of ``blocks'' and ``barriers'' will be used in the later
sections, for the definition of ``Nash-Williams orderings.'' 
For $X\in [\N]^\omega$ and $a\in\N$, we define $X^{>a}\in [\N]^\omega$
by $X^{>a}:=\{x\in X:x>a\}$.

\begin{definition}
A subset $B$ of $[X]^{<\omega}$ is called
a {\bf block} (with {\bf base $X$}) if
\begin{itemize}
\item[(B$_1$)]
$[X]^\omega = \bigcup_{s\in B}\ [s,X^{>\max s}]$, and
\item[(B$_2$)]
$s\not\propini t$ for all $s,t\in B$.
\end{itemize}
\end{definition}

In other words,
a subset $B$ of $[X]^{<\omega}$ is a block with base $X$ if and only if 
for every strictly increasing infinite sequence $(x_0,x_1,\dots)$
of elements of $X$, there exists a unique $n$ 
such that $(x_0,x_1,\dots,x_n)\in B$. 
It follows that every strictly increasing sequence over 
$X$ (finite or infinite) has at most one initial segment which lies in $B$.
(Since every element of $[X]^{<\omega}$ occurs as an initial segment of some
element of $[X]^\omega$.)
Moreover, 
note that for a block $B$, its base is given by $X=\bigcup\{b:b\in B\}$.

\begin{examples}\ 
\begin{enumerate}
\item For any $n>0$, $[X]^n$ is a block with
base $X$.
\item Suppose $B$ is a block with base $X$, and let $C\tm B$,
$Y\in [X]^\omega$.
Then $C$ is a block with base $Y$ if and only
if $C=B\cap [Y]^{<\omega}$.
\end{enumerate}
\end{examples}

A block $B$ is called a {\bf barrier} if instead of (B$_2$) it
satisfies the stronger condition
\begin{itemize}
\item[(B$_2'$)] $s\not\subset t$ for all $s,t\in B$.
\end{itemize}
The statements in the examples above remain true if ``block'' is
replaced by ``barrier''. 
Theorem~\ref{GPE}  has the following
important consequence:

\begin{cor}\label{FundamentalLemma}
If $B$ is a block with base $X$
and $B=B_1\cup B_2$, then there is a block $B'\tm B$
such that $B'\tm B_1$ or $B'\tm B_2$.
\end{cor}
\begin{proof}
We may assume that $B_1$ and $B_2$ are disjoint, so
$$P=\bigcup_{s\in B_1} [s,X^{>\max s}], \quad
             Q=\bigcup_{s\in B_2} [s,X^{>\max s}]$$
give a partition of $[X]^\omega$ into clopen sets $P$ and $Q$.  
Hence by Theorem~\ref{GPE}
there exists $H\in [X]^\omega$ with $[H]^\omega\tm P$ or
$[H]^\omega\tm Q$, and it follows that $B_1\cap [H]^{<\omega}$ or
$B_2\cap [H]^{<\omega}$, respectively, is a block. 
\end{proof}

Since any block contained in a barrier is itself a barrier, the corollary
remains true with ``block'' replaced by ``barrier.'' 

\subsection*{Construction of barriers}
The following construction turns out to be very useful:

%
%
%

\begin{prop}\label{B2}
If $B$ is a barrier with base $X\in [\N]^\omega$, then so is
$$B(2) := \bigl\{b(1)\cup b(2) : b(1),b(2)\in B, b(1)\tri b(2)\bigr\}.$$
For every $b\in B(2)$ there exist unique $b(1),b(2)\in B$ such that
$b=b(1)\cup b(2)$ and $b(1)\tri b(2)$.
\end{prop}

In the proof, we need:

\begin{lemma}\label{TriLemma}
For every $Y\in [X]^\omega$ and $b(1)\in B$ with $b(1)\ini Y$
there exists $b(2)\in B$
such that $b(1)\tri b(2)$ and $b(1)\cup b(2) \ini Y$.
\end{lemma}
\begin{proof}
The infinite sequence $Y\ohne\{\min Y\}\in [X]^\omega$
has an initial segment $b(2)$ in $B$, by (B$_1$).
By (B$_2$), $b(2) \not\propini b(1)$, and therefore 
$b(1)\tri b(2)$.
\end{proof}

\begin{proof}[Proof \rom{(of Proposition~\ref{B2})}]
Let $Y\in [X]^\omega$, and let $b(1)$ be an initial segment of $Y$ in $B$.
By the lemma, there exists $b(2)\in B$ with $b(1)\tri b(2)$
and $b(1)\cup b(2)\ini Y$. Thus $B(2)$ satisfies (B$_1$).
Now let $b(1),b(2),c(1),c(2)\in B$, with 
$b(1)\tri b(2)$, $c(1)\tri c(2)$,  and suppose that
$b\tm c$, where
$b=b(1)\cup b(2)$, $c=c(1)\cup c(2)$. 
Then
$b(2)=b\ohne\{\min b\}\tm c\ohne\{\min c\}=c(2)$, hence
$b(2)=c(2)$, by (B$_2'$). It follows that $b=c$ and $c(1)=b(1)$.
This shows that $B(2)$ is a barrier, and
also implies the last statement.
\end{proof}

\begin{cor}\label{Cstar}
Let $B$ be a barrier with base $X$. 
\begin{enumerate}
\item If $b,c\in B(2)$, then $b\tri c$ if and only if $b(2)=c(1)$. 
\item If $C\tm B(2)$ is a barrier with base $Y\tm X$, then 
$$C^*:=\bigl\{c(1),c(2):c\in C\bigr\}\tm B$$ is a barrier with base $Y\tm X$. 
If
$c_1,c_2\in C^*$ satisfy $c_1\tri c_2$, then there exists a unique $c\in C$
such that $c_1=c(1)$ and $c_2=c(2)$.
\end{enumerate}
\end{cor}
\begin{proof}
Part (1) follows from (B$_2$) and the fact that
$b(2)=b\ohne\{\min b\}$ and $c(1)\ini c\ohne\{\max c\}$.
For part (2), suppose $C\tm B(2)$ is a barrier with base $Y\tm X$, so
$C=B(2)\cap [Y]^{<\omega}$. Let $b(1)\in B\cap [Y]^{<\omega}$. By
Lemma~\ref{TriLemma}, there exists $b(2)\in B\cap [Y]^{<\omega}$ with
$b(1)\tri b(2)$. Hence $b=b(1)\cup b(2)$ is an element of 
$B(2)\cap [Y]^{<\omega}=C$, and thus $b(1),b(2)\in C^*$. This shows
that $C^*=B\cap [Y]^{<\omega}$, that is, $C^*$ is a barrier.
Finally, if $c_1\tri c_2$ are elements of $C^*$, let $c=c_1\cup c_2\in B(2)$;
then $c(1)=c_1$, $c(2)=c_2$ as required.
\end{proof}

\section{Orderings}\label{Orderings-Section}

\noindent
We state some definitions and facts concerning (partially)
ordered sets and maps between them, and give our key examples.

\subsection*{Ordered sets}

A {\bf quasi-ordering} on a set
$S$ is a binary relation $\leq$ on $S$ which is reflexive and
transitive; in this case, we call $(S,\leq)$ a {\bf quasi-ordered
set.} (If no confusion is possible, we will omit $\leq$ from the notation,
and just call $S$ a quasi-ordered set.)
If in addition the relation $\leq$ is anti-symmetric,
then $\leq$ is called an {\bf ordering} on the set $S$, and
$(S,\leq)$ (or $S$) is called an {\bf ordered set.}
If moreover $x \leq y$ or $y \leq x$ for
all $x , y \in S$, then $\leq$ is called a {\bf total ordering} on
$S$. 
We write as usual $x<y$ if $x\leq y$ and $y\not\leq x$.

\subsection*{Maps between ordered sets}
A function $\varphi\colon S\to T$ between quasi-ordered sets $(S,{\leq}_S)$ and
$(T,{\leq}_T)$ is called {\bf increasing} if $$x\leq_S y \Rightarrow
\varphi(x)\leq_T\varphi(y),\qquad\text{for all $x,y\in S$,}$$ 
and {\bf strictly increasing} 
if
$$x <_S y \Rightarrow
\varphi(x) <_T\varphi(y),\qquad\text{for all $x,y\in S$.}$$
Given a quasi-ordering $\leq$ on a set $S$, there exists a
unique ordering on the set $S/{\sim}=\{a/{\sim}:a\in S\}$ 
of equivalence classes of the equivalence
relation $$x\sim y 
\qquad\Longleftrightarrow\qquad x\leq y \text{ and } y\leq x$$ on 
$S$ such that the surjective map $a\mapsto a/{\sim}\colon S\to S/{\sim}$ is
increasing. 
Hence there is usually no loss in generality when working
with orderings rather than quasi-orderings. In the following, we shall
therefore concentrate on {\sl ordered}\/ sets, and mostly leave it to
the reader to adapt the definitions and results to the quasi-ordered case. 

\subsection*{Quasi-embeddings, embeddings, and isomorphisms}
A {\bf quasi-em\-bed\-ding}  between ordered sets $(S,{\leq}_S)$ and
$(T,{\leq}_T)$ is a map $\varphi\colon S\to T$ such that
$$\varphi(x)\leq_T \varphi(y) \Rightarrow
x\leq_S y,\qquad\text{for all $x,y\in S$,}$$
and if in addition $\varphi$
is increasing, then $\varphi$ is called an {\bf embedding.}
Any quasi-embedding between ordered sets
is injective, and any embedding is strictly increasing.
A surjective embedding $S\to T$ is an {\bf isomorphism} 
of the ordered sets $(S,{\leq}_S)$ and $(T,{\leq}_T)$.
We say that an ordering $\preceq_S$  on $S$ 
{\bf extends}
the ordering ${\leq}_S$ if the identity on $S$
is an increasing map between $(S,{\leq}_S)$ and 
$(S,{\preceq}_S)$, that is, if ${\leq}_S\subseteq {\preceq}_S$
(as subsets of $S\times S$).
We write $(S,{\leq_S})\subseteq (T,{\leq_T})$ if
$S\subseteq T$ and the natural inclusion $S\to T$ is an embedding
(i.e., $\leq_S$ is the restriction of $\leq_T$ to $S$).

\begin{examplesNumbered}\label{Ordering-Examples}
Here are some methods for constructing new ordered sets from old ones.
\begin{enumerate}
\item Any subset of an ordered set $(S,\leq)$ can be 
naturally made into an ordered set
by restricting the ordering $\leq$ to this subset.
\item The disjoint union $S \amalg T$ of ordered sets
$(S,\leq_S)$ and $(T,\leq_T)$ can be naturally made
into an ordered set via the relation $\leq_S\cup\leq_T$.
\item The cartesian product $S\times T$ of ordered sets
$(S,\leq_S)$ and $(T,\leq_T)$ can be naturally
made into an ordered set by the product ordering
$$(x,y)\leq_{S\times T} (x',y') \quad:\iff\quad x\leq_S x' \wedge
y\leq_T y'.$$
\item Let $(S,\leq)$ be an ordered set.
The free commutative monoid $S^\diamond$ generated by $S$ is naturally
ordered as follows: $s_1\cdots s_m \leq^\diamond t_1\cdots t_n :\iff$ 
there exists an
injective map $\varphi\colon \{1,\dots,m\}\to\{1,\dots,n\}$ such that
$s_i\leq_S t_{\varphi(i)}$ for all $i=1,\dots,m$.
\end{enumerate}
\end{examplesNumbered}

\subsection*{Key examples}
Two examples will play an important role in further
sections. The first one is the set 
$\N^n$ of $n$-tuples of natural
numbers, ordered by the product
ordering. It is sometimes convenient to identify the elements
of $\N^n$ with {\sl monomials}\/ in a polynomial ring as follows:
Let $X=\{X_1,\dots,X_n\}$ be a set of distinct indeterminates
and $X^\diamond=\{X^\nu: \nu\in\N^n\}$ be 
the free commutative monoid generated by $X$, where
$X^\nu:=X_1^{\nu_1}\cdots X_n^{\nu_n}$ for $\nu=(\nu_1,\dots,\nu_n)\in\N^n$.
Order $X^\diamond$ by divisibility: 
$$X^\nu\leq X^\mu \quad :\Longleftrightarrow \quad
\mu=\nu+\lambda \text{ for some $\lambda\in\N^n$.}$$ 
That is, in the context of example (4) above (for $S=X$), 
the ordering on $X^\diamond$ is $\leq^\diamond$ where $\leq$ is the
trivial ordering on $X$.
The map $\nu\mapsto X^\nu\colon\N^n\to X^\diamond$ is an
isomorphism of ordered sets. The elements of $X^\diamond$ can be seen 
as the monomials in the polynomial ring $R[X]=R[X_1,\dots,X_n]$ over an
arbitrary commutative ring $R$.

Our second example is the ordering $\red$ on the set $\Z^n$ of
$n$-tuples of integers defined as follows: if $a=(a_1,\dots,a_n)$ and
$b=(b_1,\dots,b_n)$ are elements of $\Z^n$,  we put
$a\red b$ if $a_jb_j\geq 0$ and 
$\abs{a_j}\leq \abs{b_j}$ for all $j=1,\ldots,n$. 
That is, $a\red b$ if and only if
$a$ belongs to the same orthant of $\R^n$ as $b$, and each of
its components is
not greater in absolute value than the corresponding component of $b$.
Note that $\red$ extends the product ordering on $\N^n\subseteq\Z^n$.
Moreover, we have $0\red a$ for all $a\in\Z^n$, and if $a\red b$, then
$-a\red -b$, for all $a,b\in\Z^n$. 
The ordered set $(\Z^n,{\red})$ can be identified
in a natural way with
the $n$-fold cartesian product of
$(\Z,\red)$ with itself. We have 
$(\Z^{\leq 0},\geq) \subseteq (\Z,\red)$ and
$(\Z^{\geq 0},\leq) \subseteq (\Z,\red)$, where
$\Z^{\leq 0}$ ($\Z^{\geq 0}$) denotes the set of non-positive integers
(non-negative integers, respectively). 
This somewhat roundabout way of describing $(\Z^n,\red)$ will come
in handy in the proof of Corollary~\ref{Power-Set-Cor}.

\subsection*{Initial and final segments}
An {\bf initial segment} of an ordered set $(S,\leq)$ is a
subset $I \subseteq S$ such that $$x \leq y \wedge y \in I \Rightarrow x
\in I, \qquad \text{for all $x , y \in S$.}$$ 
Dually, $F \subseteq S$ is called a
{\bf final segment} if $S \backslash F$ is an initial segment.
Given an arbitrary subset $X$ of $S$, we
denote by $$(X) := \bigl\{y \in S : \exists x \in X ( x \leq y ) \bigr\}$$
the final segment of $S$ generated by $X$, 
and by
$$[X] := \bigl\{y \in S : \exists x \in X ( x \geq y ) \bigr\}$$
the initial segment generated by $X$.
The set ${\cal I}(S)$ of initial segments of $S$ is naturally 
ordered by inclusion.
Dually, we construe the set ${\cal F}(S)$ of final segments of $S$ as
an ordered set, with the ordering given by {\sl reverse}\/ inclusion.
The intersection and union of an arbitrary family of initial (final) 
segments of $S$ is an initial (resp., final) segment of $S$.

\begin{example}
The isomorphism $\nu\mapsto X^\nu$ identifies $\N^n$ and the monoid
of monomials in the polynomial ring $R[X]$. 
The ordered set $\bigl({\cal F}(\N^n),\supseteq\bigr)$  of 
final segments of $\N^n$ may  be identified with
the set of
{\sl monomial ideals}\/ of $R[X]$ (that is, ideals of $R[X]$ which are
generated by monomials), ordered by reverse inclusion. 
\end{example}

Given a quasi-ordered set $(S,\leq)$, we can define a quasi-ordering
$\leq_{{\cal P}(S)}$ on the power set ${\cal P}(S)$ of $S$ as follows: 
for $X,Y\in {\cal P}(S)$,
\begin{equation}\label{Power-Set}
X\leq_{{\cal P}(S)} Y\quad\Longleftrightarrow\quad
\text{for every $y\in Y$ there exists $x\in X$
such that $x\leq y$.}
\end{equation} 
The map $F\mapsto F/{\sim}
\colon {\cal F}(S)\to{\cal P}(S)/{\sim}$ is an isomorphism
of ordered sets. Here $\sim$ is the equivalence relation on ${\cal P}(S)$
associated to $\leq_{{\cal P}(S)}$ as in the beginning of this section.

\subsection*{Pullback of final segments}
For any function $\varphi\colon S\to T$ between ordered sets
$(S,{\leq}_S)$ and $(T,{\leq}_T)$, we get an induced increasing function
$$\varphi^*\colon F\mapsto\bigl(\varphi^{-1}(F)\bigr),\ {\cal F}(T)\to{\cal F}(S)$$
between the ordered sets  $\bigl(\mathcal{F} ( T ),\supseteq\bigr)$
and  $\bigl(\mathcal{F} ( S ),\supseteq\bigr)$.

\begin{example}
Let $(S,{\leq}_S)$,
$(T,{\leq}_T)$, and $(U,{\leq_U})$ be ordered sets, and
suppose that $(S,{\leq}_S),(T,{\leq_T})\subseteq (U,{\leq_U})$, so the
natural inclusions
$i_S\colon S\to U$ and $i_T\colon T\to U$ are increasing.
Then $$i_S^*(F)=F\cap S, \quad
i_T^*(F)=F\cap T\qquad\text{for any $F\in{\cal F}(U)$.}$$ 
Hence
if in addition $U=S\cup T$, then
$i_S^*\times i_T^*$ gives an embedding
${\cal F}(U) \to
{\cal F}(S)\times {\cal F}(T)$.
\end{example}

The following rules will be used later:

\begin{lemma}\label{FRules}
Let $\varphi\colon S\to T$.
\begin{enumerate}
\item If  $\varphi$ is a quasi-embedding, 
then $\varphi^*$ is surjective.
\item If $\varphi$ is increasing and surjective, then
$\varphi^*$ is a quasi-embedding.
\end{enumerate}
\end{lemma}
\begin{proof}
Since part (2) is clear, we just prove (1).
Suppose $\varphi$ is a quasi-embedding, and
let $G\in{\cal F}(S)$. Let $F$ be the final segment of
$T$ generated by $\varphi(G)$.  Then clearly $G\tm\varphi^{-1}(F)$.
Conversely, if $x\in S$ satisfies $\varphi(x)\in F$, then $\varphi(x)\geq
\varphi(g)$ for some $g\in G$. Since $\varphi$ is a quasi-embedding,
we have $x\geq g$, so $x\in G$. This shows $\varphi^*(F)=G$. 
\end{proof}

\subsection*{Antichains}
For elements $x,y$ of an ordered set $S$, we write
$x \| y$ if $x\not\leq y$ and $y\not\leq x$.
An {\bf antichain} of $S$ is a subset $A \subseteq S$
such that $x \| y$ for all $x \neq y$ in $A$. 
An element $x$ of $S$ is called a {\bf minimal} element of $S$
if $y\leq x \Rightarrow y=x$
for all $y\in S$. The minimal elements of a subset $X\tm S$ 
form an antichain, denoted by $X_{\min}$. 

\subsection*{Well-founded orderings}
An ordered set $S$ is {\bf well-founded} if there is no
infinite strictly decreasing sequence $x_0 > x_1 > \cdots $ in $S$. 
For any element $x$ of a subset $X$ 
of a well-founded ordered set $S$, 
there is at least one minimal element $y\in X_{\min}$ with $y\leq x$. 
It follows that any final segment $F$ of $S$ is generated by its antichain
$F_{\min}$ of minimal elements; in fact, $F_{\min}$ is the (unique)
smallest generating set for $F$.

\section{Noetherian Orderings}

\noindent
We say that an
ordered set $S$ is {\bf Noetherian} if it is well-founded
and every antichain of $S$ is finite.
More generally, 
we say that a quasi-ordered set $S$ is Noetherian if the associated
ordered set $S/{\sim}$ is Noetherian. 
Since every antichain of a
totally ordered set consists of at most one
element, a totally ordered set $S$ is Noetherian if and only if it is
well-founded; in this case $S$ is called {\bf well-ordered.}

\begin{remark}
Noetherian orderings are usually called ``well-quasi-orderings'' in the
literature (see, e.g., \cite{Kruskal}). Following a proposal by
Joris van der Hoeven \cite{vdH:phd} we use the more concise
term ``Noetherian''. 
\end{remark}

An infinite sequence 
$x_0,x_{1}, \ldots$ in $S$ is {\bf good} if $x_i \leq x_j$ for some $i <
 j$, and {\bf bad,} otherwise. Clearly, if $\{x_0,x_1,\dots\}$ is an
antichain, then $x_0,x_1,\dots$ is a bad sequence.
The following characterization of Noetherian orderings is folklore 
(see, e.g., \cite{Milner}).

\begin{prop}\label{Folk}
Let $S$ be an ordered set. The following are equivalent:
\begin{enumerate}
\item $S$ is Noetherian.
\item Every infinite sequence $x_0 , x_1 , \ldots$ in $S$ contains an
increasing subsequence.
\item Every infinite sequence $x_0 , x_1 , \ldots$ in $S$ is good.
\item Any subset $X\tm S$ has only finitely many minimal elements, and
for every $x\in X$ there is a minimal element $y$ of $X$ with $x\geq y$. 
\item Any final segment of $S$ is finitely generated.
\item $\bigl(\mathcal{F} ( S ),\supseteq\bigr)$ is well-founded \textup{(}i.e., the ascending chain
condition with respect to inclusion holds for final segments of $S$\textup{).}
\item $\bigl(\mathcal{I} ( S ),\subseteq\bigr)$ is well-founded \textup{(}i.e., the descending chain
condition with respect to inclusion holds for initial segments of $S$\textup{).}\end{enumerate}
\end{prop}
\begin{proof}
The implication (1)~$\Rightarrow$~(2) follows from applying Ramsey's
Theorem to the partition $[\N]^2=P\cup Q\cup R$, where
$$P = \bigl\{ \{i,j\} : x_i\| x_j\bigr\}, \quad
Q = \bigl\{ \{i,j\} : i<j, x_i > x_j\bigr\},$$
and $R=[\N]^2\ohne (P\cup Q)$.
The implications (2)~$\Rightarrow$~(3)~$\Rightarrow$~(1)
are trivial, and 
(1)~$\Rightarrow$~(4)~$\Rightarrow$~(5) follows from the remarks at the end
of the last section.
If $F_0\tm F_1\tm \cdots$ is an ascending chain of final segments of $S$,
then $F=\bigcup_n F_n$ is a final segment of $S$. If $F$ is
finitely generated, say by $X\tm F$, then $X\tm F_n$ for some $n$;
thus $F_n=F_{n+1}=\cdots$. This shows (5)~$\Rightarrow$~(6);
by passing to complements, we obtain (6)~$\iff$~(7). For
(6)~$\Rightarrow$~(3), let $x_0,x_1,\dots$ be a sequence in $S$. 
By (6), the sequence $(x_0)\tm (x_0,x_1)\tm \cdots$
of final segments of $S$ becomes stationary: for some $n$, we have
$x_j\in (x_0,\dots,x_n)$ for all $j>n$. In particular,
$x_i\leq x_{n+1}$ for some
$i\in\{0,\dots,n\}$.
\end{proof}

The proposition now immediately provides the following construction methods
for Noetherian orderings:

\begin{examples}
Suppose $(S,{\leq}_S)$ and
$(T,{\leq}_T)$ are ordered sets. Then:
\begin{enumerate}
\item If there exists an increasing surjection $S\to T$, and
$(S,{\leq}_S)$ is Noetherian, then $(T,{\leq}_T)$ is Noetherian.
In particular, if $(S,{\leq}_S)$ is Noetherian, then 
any ordering on $S$ which extends $\leq_S$ is
Noetherian.
\item If there exists a quasi-embedding $S\to T$, and $(T,{\leq}_T)$
is Noetherian, then $(S,{\leq}_S)$ is Noetherian.
In particular, if $(T,{\leq}_T)$ is Noetherian, then any subset of $T$ with 
the induced ordering is Noetherian.
\item  If $(S,{\leq}_S)$ and $(T,{\leq}_T)$ are Noetherian and
$(U,{\leq}_U)$ is an ordered set with
 $(S,{\leq}_S),(T,{\leq_T})\subseteq (U,{\leq_U})$, then
$S\cup T$ is Noetherian. 
In particular, it follows that  $S\amalg T$ is Noetherian.
\item If $(S,{\leq}_S)$ and $(T,{\leq}_T)$ are Noetherian, then so is
$S\times T$.
Inductively, it follows that 
if the ordered set $(S,{\leq}_S)$ is Noetherian, 
then so is $S^n$ equipped with the product ordering, for
every $n>0$.
\end{enumerate}
\end{examples}

Applying Example~(4) to $S=\N$, we obtain:

\begin{cor}
\rom{(Dickson's Lemma.)} For each $n>0$, the ordered set $\N^n$ is
Noetherian. \qed
\end{cor}

By combining Examples~(3) and (4) with the analysis of the ordered
set $(\Z^n,{\red})$ given in the last section, we get:

\begin{cor}\label{red-Noetherian}
For each $n>0$, the ordered set $(\Z^n,\red)$ is Noetherian. \qed
\end{cor}

\section{Strongly Noetherian Orderings}

\noindent
Besides the results stated in the
examples following Proposition~\ref{Folk}, several other preservation
theorems for Noetherian orderings are known. For example, if $(S,\leq)$ is a
Noetherian ordered set, then so is $(S^\diamond,\leq^\diamond)$ as
defined in Section~\ref{Orderings-Section}. 
(This was first 
proved by Higman \cite{Higman}, with a simplified proof given by
Nash-Williams \cite{NW2}.) Similar theorems can be proved for
other ordered sets built from $S$, for example the collection of 
all finite trees whose nodes are labelled by elements of $S$, ordered by
the homeomorphic embedding ordering (Kruskal \cite{Kruskal2}).
The common feature of all these constructions
is their {\it finitary}\/ character.
If one builds ordered sets by allowing operations
of an infinite nature, the situation changes drastically:
for example, if $S$ is a Noetherian ordered set, 
then $\bigl({\cal F}(S),\supseteq\bigr)$
is in general {\it not}\/ Noetherian. An example for this phenomenon was
first given by Rado \cite{Rado}: Let
$R:=[\N]^2$,
ordered by the rule
$$(i,j)\leq_R (k,l) \qquad\iff\qquad
\text{either $i=k$ and $j\leq l$, or $j<k$.}$$
It is quickly verified (see, e.g., \cite{Maclagan}) 
that $(R,\leq_R)$ is a Noetherian ordered set, but
$\bigl({\cal F}(R),\supseteq\bigr)$ is not Noetherian: the sequence
$F_1,F_2,\dots$ where $F_j$ is the final segment of $R$ generated by
all $(i,j)\in\N^2$ with $i<j$, is an infinite antichain in ${\cal F}(R)$.
This example is archetypical in the following sense:

\begin{theorem}\label{Rado}
For a Noetherian ordered set $(S,\leq_S)$, the following are equivalent:
\begin{enumerate}
\item $\bigl({\cal F}(S),\supseteq\bigr)$ is not Noetherian.
\item $\bigl({\cal I}(S),\subseteq\bigr)$ is not Noetherian.
\item There exists a function $f\colon [\N]^2\to S$ such that
$f(i,j)\not\leq_S f(j,k)$ for all $i<j<k$.
\item There exists an embedding $(R,\leq_R)\to (S,\leq_S)$.
\item There exists a quasi-embedding $(R,\leq_R)\to (S,\leq_S)$.
\end{enumerate}
\end{theorem}
\begin{proof}
The implication (1)~$\Rightarrow$~(2)
is trivial. Suppose $I_0,I_1,\dots$ is a bad sequence in  
$\bigl({\cal I}(S),\subseteq\bigr)$. So for each $i<j$ there
exists $a\in I_i$ with $a\not\leq b$ for all $b\in I_j$. Hence it is
possible to choose, for each $i<j$, an element $f(i,j)\in I_i$ such
that for all $i<j<k$, we have $f(i,j)\not\leq_S f(j,k)$. 
This shows 
(2)~$\Rightarrow$~(3). For (3)~$\Rightarrow$~(4), let
$f\colon [\N]^2\to S$, $f(i,j)=a_{ij}$, be as in (3).
Consider
the partitions
$$[\N]^3=P\cup Q, \qquad [\N]^4=P'\cup Q'$$
given by
$$P=\bigl\{\{i,j,k\}: i<j<k, a_{ij}\leq a_{ik}\bigr\}, \qquad Q=[\N]^3\ohne P$$
and
$$P'=\bigl\{\{i,j,k,l\}: i<j<k<l, a_{ij}\leq a_{kl}\bigr\}, \qquad Q'=[\N]^4\ohne P'.$$
By applying Ramsey's Theorem twice we obtain an infinite set $H\tm\N$
which is homogeneous for both partitions. Since $S$ is Noetherian, we must
have $[H]^3\tm P$, $[H]^4\tm P'$. It follows that 
$$(i,j)\leq_R (k,l) \qquad\iff\qquad a_{ij}\leq_S a_{kl},$$ 
for all $i<j$, $k<l$ in $H$. Therefore, $\{a_{ij}:i<j, i,j\in H\}$
(with the ordering induced from $S$) is isomorphic to $(R,\leq_R)$. 
The implication (4)~$\Rightarrow$~(5) is again trivial, and
(5)~$\Rightarrow$~(1) follows from Lemma~\ref{FRules},~(1) and
Example (1) following Proposition~\ref{Folk}.
\end{proof}

Let us call an ordered set $(S,\leq_S)$ {\bf strongly Noetherian} if
$\bigl({\cal F}(S),\supseteq\bigr)$ is Noetherian.
(So {\it a forteriori,}\/ $(S,\leq_S)$ is Noetherian, by 
Proposition~\ref{Folk}.)
The following facts follow easily from Lemma~\ref{FRules} and the
examples after Proposition~\ref{Folk}:

\begin{examples}
Suppose $(S,{\leq}_S)$ and
$(T,{\leq}_T)$ are ordered sets.
Then:
\begin{enumerate}
\item If there exists an increasing surjection $S\to T$ and
$(S,{\leq}_S)$ is strongly Noetherian, then $(T,{\leq}_T)$ is strongly
Noetherian.
In particular, if $(S,{\leq}_S)$ is strongly Noetherian, then 
any ordering on $S$ which extends $\leq_S$ is
strongly Noetherian.
\item If there exists a quasi-embedding $S\to T$ and $(T,{\leq}_T)$
is strongly Noetherian, then $(S,{\leq}_S)$ is strongly Noetherian.
In particular, if $(S,{\leq}_S)$ is strongly Noetherian, 
then any subset of $S$ with the induced ordering is strongly Noetherian.
\item If $(S,{\leq}_S)$ and $(T,{\leq}_T)$ are strongly Noetherian and
$(U,{\leq}_U)$ is an ordered set with
 $(S,{\leq}_S),(T,{\leq_T})\subseteq (U,{\leq_U})$, then
$S\cup T$ is strongly Noetherian. 
In particular, it follows that  
if $(S,{\leq}_S)$ and $(T,{\leq}_T)$ are strongly 
Noetherian, then so is $S\amalg T$.
\end{enumerate}
\end{examples}

Strong Noetherianity is also preserved under cartesian products:

\begin{prop}\label{CP}
If $(S,{\leq}_S)$ and $(T,{\leq}_T)$ are strongly 
Noetherian ordered sets, then $S\times T$ is strongly Noetherian.
\end{prop}

(This fact was stated without proof in \cite{Maclagan} and attributed
there to Farley and Schmidt. For another proof see \cite{AschenbrennerPong}.)

\begin{proof}
Let $f\colon [\N]^2\to S\times T$.
We denote by  $\pi_S\colon S\times T\to S$ the projection $(s,t)\mapsto s$
onto the
first component, and we put $f_S:=\pi_S\circ f$. Similarly we define
$f_T\colon [\N]^2\to T$. 
 Since $S$ is strongly
Noetherian, there are $i<j<k$ with $f_S(i,j)\leq_S f_S(j,k)$.
(By the previous theorem.)
Now Consider the partition $[\N]^3=P\cup Q$, where
$$P=\bigl\{ \{i,j,k\}: i<j<k,\ f_S(i,j) \leq_S f_S(j,k)\bigr\}, 
\qquad Q=[\N]^3\ohne P.$$
By Ramsey's Theorem, we find an infinite homogeneous set $H\tm\N$ 
for this partition.
Since $S$ is strongly Noetherian, we must have $[H]^3\tm P$.
Changing from $\N$ to $H$, we thus may assume that
$f_S(i,j)\leq_S f_S(j,k)$ {\it for all $i<j<k$.}\/ Since $T$ is strongly
Noetherian, there are $i<j<k$ with $f_T(i,j)\leq_T f_T(j,k)$. Hence
$f(i,j)\leq f(j,k)$. Thus $S\times T$ is strongly Noetherian.
\end{proof}

\begin{cor}\label{Maclagan} \rom{(Maclagan, \cite{Maclagan}.)}
The ordered set $\N^n$ is strongly Noetherian, for every $n>0$. \qed
\end{cor}


\section{Nash-Williams Orderings}\label{Section-Better-Noetherian}

\noindent
By the equivalence of (1) and (3) in Theorem~\ref{Rado}, 
strong Noetherianity may be expressed using the
terminology introduced in Section~\ref{Section-Preliminaries}:
{\it An ordered set
$S$ is strongly Noetherian if and only if for every 
function $f\colon B\to S$, where $B\tm [\N]^2$ is a barrier, 
there exist $b_1,b_2\in B$
with $b_1\tri b_2$ and $f(b_1)\leq f(b_2)$.}\/
The search for a combinatorial condition on an ordered set
$S$ which ensures that  
$\bigl({\cal F}(S),\supseteq\bigr)$ is not only Noetherian, but
{\it strongly}\/ Noetherian, therefore naturally leads to the following
concept introduced by Nash-Williams \cite{NW}  (under the name of 
``better well-quasi-ordering'').
Below, we shall call a function $f\colon B\to S$, whose domain $B\tm [\N]^{<\omega}$ 
is a barrier, an {\bf $S$-array.} (So in particular, every sequence
$x_0,x_1,\dots$ of elements of $S$ can be considered as an $S$-array.)
We say that an $S$-array $f\colon B\to S$ is {\bf good} if there are
$b_1,b_2\in B$ such that $b_1\tri b_2$ and $f(b_1)\leq f(b_2)$, 
{\bf bad} if it is not good,
and {\bf perfect} if $f(b_1)\leq f(b_2)$ {\it for all}\/ $b_1\tri b_2$ in $B$.

\begin{definition}
An ordered set $S$ is {\bf Nash-Williams} if every $S$-array is good.
A quasi-ordered set $S$ is Nash-Williams if 
$S/{\sim}$ is a Nash-Williams ordered set.
\end{definition}

%
Clearly if $S$ is Nash-Williams, then $S$ is Noetherian. Moreover:

\begin{lemma}
If $S$ is well-ordered, then $S$ is Nash-Williams.
\end{lemma}
\begin{proof}
Suppose $f\colon B\to S$ is 
an $S$-array. Let $b_1\in B$ be such that
$f(b_1)=\min f(B)$. By Lemma~\ref{TriLemma} there exists $b_2\in B$
such that $b_1\tri b_2$, and we have $f(b_1)\leq f(b_2)$.
\end{proof}

As we did for strong Noetherianity, we will now successively show that each of
the constructions exhibited in the examples following Proposition~\ref{Folk}
preserves the Nash-Williams property as well.

\begin{lemma}
Let $(S,{\leq_S})$ and $(T,{\leq_T})$ be ordered sets.
\begin{enumerate}
\item If there exists an increasing surjection $S\to T$ and
$S$ is Nash-Williams, then $T$ is Nash-Williams.
\item If there exists a quasi-embedding $S\to T$ and
$T$ is Nash-Williams, then $S$ is Nash-Williams.
\end{enumerate}
\end{lemma}
\begin{proof}
Let $\varphi\colon S\to T$.
For part (1), suppose $\varphi$ is increasing and surjective, and
let $f\colon B\to T$ be
a $T$-array. Choose any function $\psi\colon T\to S$ such that
$\varphi\circ\psi=\id_T$. Then $\psi\circ f\colon B\to S$ is an $S$-array.
Since $S$ is Nash-Williams, there exist $b_1,b_2\in B$ with
$b_1\tri b_2$ and $\psi(f(b_1))\leq_S\psi(f(b_2))$, hence $f(b_1)\leq_S
f(b_2)$. For part (2), assume that $\varphi$ is a quasi-embedding, and
let $g\colon B\to S$ be an $S$-array. Then
$\varphi\circ g\colon B\to T$ is a $T$-array. 
Since $T$ is Nash-Williams, there exist $b_1,b_2\in B$ with
$b_1\tri b_2$ and $\varphi(g(b_1))\leq_T\varphi(g(b_2))$, hence $g(b_1)\leq_S
g(b_2)$, since $\varphi$ is a quasi-embedding.
\end{proof}

In particular, if $(S,{\leq})$ is a Nash-Williams ordered set, 
then so is any subset of $S$ with the induced ordering, and
any ordering on $S$ which extends the ordering $\leq$ on $S$.

\begin{lemma} \label{Union-Better-Noetherian}
Suppose that $(S,{\leq}_{S})$, $(T,{\leq}_{T})$ and $(U,\leq_U)$ are 
ordered sets, with $(S,{\leq}_{S}),(T,{\leq}_S)\tm (U,\leq_U)$. If
$S$ and $T$ are Nash-Williams, then
$S\cup T$ is also Nash-Williams.
\end{lemma}
\begin{proof}
Let $f\colon B\to S\cup T$ be a bad $S\cup T$-array. Let
$B_S:=f^{-1}(S)$ and $B_T:=f^{-1}(T)$. Then there exists a barrier $B'\tm B$
such that $B'\tm B_S$ or $B'\tm B_T$, by Corollary~\ref{FundamentalLemma}. 
So either $f|B'$
is a bad $S$-array, or $f|B'$ is a bad $T$-array, 
a contradiction.
\end{proof}

It follows that if $(S,{\leq}_S)$ and $(T,{\leq}_T)$
are Nash-Williams, then 
$S\amalg T$ is Nash-Williams.
Every finite ordered set is Nash-Williams.

The following fact distinguishes Nash-Williams
orderings (among Noetherian, or strongly Noetherian, orderings).

\begin{prop}\label{FinIni}
Suppose that the ordered set $S$ is  Nash-Williams.
Then the ordered sets $\bigl({\cal F}(S),\supseteq\bigr)$ 
and $\bigl({\cal I}(S),\subseteq\bigr)$ are Nash-Williams.
\end{prop}
\begin{proof}
Let $f\colon B\to {\cal I}(S)$ be a bad ${\cal I}(S)$-array. 
By Proposition~\ref{B2}, $B(2)$ is a barrier.
We construct
a bad $S$-array $g\colon B(2)\to S$ as follows: 
For every $b\in B(2)$, we have
$f(b(1)) \not\subseteq f(b(2))$; so we can 
choose $g(b)\in f(b(1))\ohne f(b(2))$.
Now suppose $c\in B(2)$ with $b\tri c$. Then $c(1)=b(2)$
(by Corollary~\ref{Cstar},~(1)), hence
$g(b)\not\leq g(c)$. Therefore, $g$ is bad. This shows that
${\cal I}(S)$ is Nash-Williams.
Suppose $h\colon C\to{\cal F}(S)$ is an ${\cal F}(S)$-array. We then
consider the ${\cal I}(S)$-array $h'\colon C\to{\cal I}(S)$ defined by
$h'(c)=S\ohne h(c)$, for all $c\in C$. Since ${\cal I}(S)$ is 
Nash-Williams, we find $c_1\tri c_2$ in $C$ with $h'(c_1)\tm h'(c_2)$, and
hence $h(c_1)\om h(c_2)$. Thus $h$ is good, and ${\cal F}(S)$ is Nash-Williams.
\end{proof}

For the definition of the quasi-ordering
$\leq_{{\cal P}(S)}$ in the statement of the next corollary,
see \eqref{Power-Set}.

\begin{cor}\label{FinIni-Cor}
If $S$ is a Nash-Williams quasi-ordered set, then
the quasi-ordered set
$\bigl({\cal P}(S),\leq_{{\cal P}(S)}\bigr)$ is  
Nash-Williams.
\end{cor}
\begin{proof}
The canonical increasing surjection $\varphi\colon S\to S/{\sim}$
induces an isomorphism $\varphi^*\colon{\cal F}(S/{\sim})\to{\cal F}(S)$.
By the proposition, ${\cal F}(S/{\sim})$, and hence ${\cal F}(S)$,
is Nash-Williams. Since ${\cal P}(S)/{\sim}\cong {\cal F}(S)$,
it follows that ${\cal P}(S)$ is Nash-Williams.
\end{proof}

Next, we want to show that the cartesian
product of any two Nash-Williams ordered sets
is again Nash-Williams. The proof is similar to the one of
Proposition~\ref{CP} and uses the following lemma.
Here, an {\bf $S$-subarray} of an $S$-array $f\colon B\to S$ is
an $S$-array $g\colon C\to S$ where $C\tm B$, $g=f|C$.

\begin{lemma}\label{BadOrPerfect}
Let $(S,\leq)$ be an ordered set.
Every $S$-array contains 
either a bad  $S$-subarray or a perfect $S$-subarray.
\end{lemma}
\begin{proof}
Let $f\colon B\to S$ be an $S$-array. Consider the 
partition $B(2)=P\cup Q$ of the barrier $B(2)$, where
$$P=\bigl\{b\in B(2): f\bigl(b(1)\bigr)\leq f\bigl(b(2)\bigr)\bigr\}, \qquad Q=B(2)\ohne P.$$
By Corollary~\ref{FundamentalLemma}, there is a barrier $C\tm B(2)$ such
that $C\tm P$ or $C\tm Q$. By Corollary~\ref{Cstar},~(2),
$C^*\tm B$ is a barrier, and $f|C^*$ is either perfect or bad, depending on
whether $C\tm P$ or $C\tm Q$.
\end{proof}

\begin{prop} \label{Prod-Better-Noetherian}
If $(S,{\leq}_S)$ and $(T,{\leq}_T)$
are Nash-Williams ordered sets, then $S\times T$ is Nash-Williams.
\rom{(\cite{NW}, Corollary 22A.)}
\end{prop}
\begin{proof}
Let $f\colon B\to S\times T$ be an $S\times T$-array; we want to show
that $f$ is good. By the previous
lemma, and since $S$ is Nash-Williams, 
the $S$-array $f_S=\pi_S\circ f$ has a perfect $S$-subarray
$B'\to S$. Restricting
$f$ to the barrier $B'$, if necessary, we may assume that $f_S\colon B\to S$ 
is perfect. Since $T$ is Nash-Williams, there exist $b_1,b_2\in B$ with
$b_1\tri b_2$ and $f_T(b_1)\leq_T f_T(b_2)$, whence $f(b_1)\leq f(b_2)$.
So $f$ is good.
\end{proof}

It follows that $\N^n$ is Nash-Williams, for any $n>0$.
Defining inductively $${\cal F}_0(\N^n):=\N^n, \qquad {\cal F}_{k+1}(\N^n):=
{\cal F}\bigl({\cal F}_k(\N^n)\bigr)\text{ for all $k\in\N$,}$$
we get, by Proposition~\ref{FinIni}:

\begin{theorem}\label{NmBetterNoeth}
Each of the sets 
$${\cal F}_1(\N^n)={\cal F}(\N^n),\ 
{\cal F}_2(\N^n)={\cal F}\bigl({\cal F}_1(\N^n)\bigr),\ \dots,$$
ordered by reverse inclusion, is Noetherian. \qed
\end{theorem}

\begin{remark}
Theorem~\ref{NmBetterNoeth} above implies that every
subset of $\cal F(\N^n)$ is strongly Noetherian, and hence yields
the theorem in the introduction.
\end{remark}

Lemma~\ref{Union-Better-Noetherian} and 
Proposition~\ref{Prod-Better-Noetherian} show that $(\Z^n,\red)$ is
Nash-Williams, for any $n>0$. As above, this implies that
each set
$${\cal F}(\Z^n,\red),\ {\cal F}\bigl({\cal F}(\Z^n,\red)\bigr),\ \dots,$$ 
ordered by reverse inclusion, is Noetherian.
Here is a slight reformulation of this fact which will be used in
Part~2.
Let $\red_0$ denote the ordering $\red$ on ${\cal P}_0(\Z^n):=\Z^n$, and
inductively 
define ${\cal P}_{k+1}(\Z^n):={\cal P}\bigl({\cal P}_k(\Z^n)\bigr)$ 
and a quasi-ordering
$\red_{k+1}$ of ${\cal P}_{k+1}(\Z^n)$ as follows:
for $X,Y\in {\cal P}_{k+1}(\Z^n)$ put
\begin{equation}\label{red_k}
X\red_{k+1} Y\quad\Longleftrightarrow\quad
\text{for each $y\in Y$ there is some
$x\in X$ with $x\red_k y$.}
\end{equation} (That is,
${\red_{k+1}}={\leq_{{\cal P}({\cal P}_{k}(\Z^n),\red_k)}}$ for all $k$,
in the notation of \eqref{Power-Set}.)
By induction on $k$, Corollary~\ref{FinIni-Cor} implies:

\begin{cor}\label{Power-Set-Cor}
The quasi-ordered 
set $\bigl({\cal P}_k(\Z^n),\red_k\bigr)$ is Noetherian, for each
$k\in\N$ and $n>0$. \qed
\end{cor}

\section{Applications to Hilbert Functions}\label{Applications to Hilbert Functions}

\noindent
As an illustration, let us show how Nash-Williams' theory
as outlined above can be employed to deduce some finiteness properties for
Hilbert functions. 
We work in a rather general setting.
Let $(S,\leq)$ be a non-empty ordered set and
$\delta\colon S\to A$ a map, where $A$ is any set. 
We think of $\delta$ as providing a
{\bf grading} of $S$ by $A$, and we call $\delta(s)\in A$ the
{\bf degree} of $s\in S$.
We call a final segment $F$ of $S$ {\bf admissible} if for any degree
$a\in A$, there are only finitely many $x\notin F$ with $\delta(x)=s$.
If the map $\delta$ has finite fibers, then every final
segment of $S$ is admissible. (This is the case, for example, if
$S$ is Noetherian and $\delta\colon S\to A$ is a
strictly increasing map into an {\sl ordered}\/ set $A$, see 
Example (1) following Proposition~\ref{Folk}.)

We extend the ordering of $\N$ to a total ordering of 
the set $\N_\infty=\N\cup\{\infty\}$ by declaring $n<\infty$ for all $n\in\N$.
For a 
final segment $F$ of $S$ and a degree $a\in A$, we let $h_F(a)$ be the
number of elements of $S$ of degree $a$ which are not in $F$ (an
element of $\N_\infty$).  We shall call
$h_F\colon A\to\N_\infty$ the {\bf Hilbert function} of 
$F\in{\cal F}(S)$. (So $F$ is admissible precisely if $h_F(a)\in\N$ for all 
$a$.) 

\subsection*{Growth of Hilbert functions}
First we show that the growth of 
Hilbert functions
is finitely determined in the following sense.
We fix an ordered set $S$ and a grading $\delta$ of $S$ by $A$ as above.

\begin{prop}\label{KeyProp}
Suppose $S$ is strongly Noetherian, that is, ${\cal F}(S)$ is Noetherian,
and let $h\colon A\to\N_\infty$ be any function. 
There exist finite sets
 ${\cal M}_{\geq h}$, ${\cal M}_{\not\leq h}$ and  ${\cal M}_{\succeq h}$ of
final segments of $S$ such
that for every final segment $F$ of $S$, we have
\begin{align*}
h_F(a) \geq h(a) \text{\ for all $a$} &\ \Longleftrightarrow\ \text{$F\subseteq E$ for
some $E\in{\cal M}_{\geq h}$}, \\
h_F(a) > h(a) \text{\ for some $a$} &\ \Longleftrightarrow\ \text{$F\subseteq E$ for
some $E\in{\cal M}_{\not\leq h}$},
\end{align*}
and
$$h_F(a) \geq h(a) \text{\ for all but finitely many $a$}\ \Longleftrightarrow\ \text{$F\subseteq E$ for
some $E\in{\cal M}_{\succeq h}$.}$$
\end{prop}
\begin{proof}
For $a\in A$, consider the subsets
$${\cal F}_a:=\bigl\{ F\in{\cal F}(S) : h_{F}(a) \geq h(a)\bigr\}, \quad
{\cal F}'_a:=\bigl\{ F\in{\cal F}(S) : h_{F}(a) > h(a)\bigr\}$$
of ${\cal F}(S)$. Since for all $E,F\in{\cal F}(S)$,
$$E\supseteq F \quad\Rightarrow\quad h_{E}(a)\leq h_{F}(a)
\text{ for all $a\in A$,}$$ 
${\cal F}_a$  and ${\cal F}_a'$ are final segments of the Noetherian
ordered set $\bigl({\cal F}(S),\supseteq\bigr)$. 
Hence there exist finite sets of generators ${\cal M}_{\geq h}$,
${\cal M}_{\not\leq h}$ and ${\cal M}_{\succeq h}$  
 for
the final segments $$\bigcap_{a} {\cal F}_a,\quad \bigcup_{a} {\cal F}'_a,\quad
\bigcup_{\substack{D\subseteq A \\ \text{finite}}}\bigcap_{a\notin D} {\cal F}_a$$ of ${\cal F}(S)$, respectively.
These sets have the required properties.
\end{proof}

The previous proposition applies to $S=\N^n$, by Corollary~\ref{Maclagan}.

\subsection*{A proposition of Haiman and Sturmfels}
Next we show how a combinatorial statement due to Haiman and Sturmfels
(\cite{Haiman-Sturmfels}, Proposition~3.2) can be obtained 
using the techniques above. This statement is crucial in 
the construction of the multigraded Hilbert scheme given in 
\cite{Haiman-Sturmfels}. 
If $S$ is strongly Noetherian, then, given any 
function $h\colon A\to\N_\infty$,
there are only finitely many admissible
final segments with
Hilbert function $h$. Moreover:

\begin{lemma}\label{HilbertDown}
Suppose $S$ is strongly Noetherian. 
Then, for any  function $h\colon A\to\N_\infty$,
there exists a finite set $D\tm A$ such that every
$F\in{\cal F}(S)$ satisfies: if $h_F(a)\leq h(a)$ for all $a\in D$, 
then $h_F(a)\leq h(a)$ for all $a\in A$.
\end{lemma}
\begin{proof}
Let $h\colon A\to\N_\infty$, and let ${\cal M}_{\not\leq h}=\{F_1,\dots,F_m\}$
be as in Proposition~\ref{KeyProp}.
For each $i$ there exists $a_i\in A$ such that $h_{F_i}(a_i)>h(a_i)$.
Let $D=\{a_1,\dots,a_m\}$, a finite set of degrees. Now suppose $F\in
{\cal F}(S)$ satisfies $h_F|D \leq h|D$. Then
$F\notin{\cal F}$: Otherwise, $F\tm F_i$ for some $i$, so $h_F(a)\geq
h_{F_i}(a)$ for all $a\in A$; but $h_{F_i}(a_i)>h(a_i)\geq h_F(a_i)$:
a contradiction. Thus $h_F(a)\leq h(a)$ for all $a$.
\end{proof}
\begin{remark}
In an analogous way, one shows:
there exists a finite set $D\subseteq A$ such that for every $F\in{\cal F}(S)$,
if $h_F(a)<h(a)$ for all $a\in D$, then $h_F(a)<h(a)$ for all $a\in A$.
In particular, there exists a finite set $D_{\text{a}}\subseteq A$ such that
for any $F\in{\cal F}(S)$,
if $h_F(a)<\infty$ for all $a\in D_{\text{a}}$, then $F$ is admissible.
\end{remark}

Note that the inequality $h_F(a)\geq h(a)$ 
is clearly not determined by finitely many degrees $a$.
(Consider $S=\N$ equipped with the grading $\delta(m)=m$, and 
$h(a)=1$ for all $a$.) 
However, we have:

\begin{lemma}\label{HilbertUp}
Suppose that $S$ is strongly Noetherian. For
any  function $h\colon A\to\N$,
there exists a finite set $D\tm A$ such that for every
final segment $F\tm S$ generated in degrees $D$: if
$h_F(a)= h(a)$ for all $a\in D$, then $h_F(a)= h(a)$ for all $a\in A$.
\end{lemma}
\begin{proof}
By induction on $i$, we
construct an increasing
sequence $(D_i)_{i\geq 0}$ of finite subsets of $A$ such that
$[C_{i+1}]\subseteq [C_i]$ for all $i$, where
$C_i$ denotes the finite antichain
consisting of all $F\in{\cal F}(S)$ which are generated in $D_i$ and
whose Hilbert function agrees with $h$ on $D_i$.
Let $C=\{F:h_F=h\}$, a finite antichain, and let
$D_0$ be a finite set of degrees such that every $F\in C$ is generated
in $D_0$.
Inductively, suppose that $D_i$ has been constructed already.
Let $D_{i+1}\supseteq D_i$ be a finite set of degrees such that
every $F\in C_i$ with $h_F|D_i=h|D_i$ is in $C_{i+1}$.
We have $[C_{i+1}]\subseteq [C_i]$: if $F\in C_{i+1}$, then the final
segment $G$ generated by the elements of $F$ 
with degrees in $D_i$ is in $C_i$, and $F\subseteq G$.
By Proposition~\ref{Folk}, $C_{i+1}=C_i$ for some $i$; $D=D_i$ works.
\end{proof}

From now on we restrict to the special case $S=\N^n$, $n>0$.
Identifying monomial ideals in the polynomial ring 
$R[X]=R[X_1,\dots,X_n]$ (where $R$ is  
a fixed commutative ring)
with final segments of $\N^n$ as usual, we obtain
the following strengthening of Proposition~3.2 in \cite{Haiman-Sturmfels}.
In that paper, $A$ is an abelian group and $\delta=\deg
\colon\N^n\to A$ a semigroup homomorphism. (So a monomial ideal $I$ of $R[X]$
is admissible exactly if $h_I(0)<\infty$.)

\begin{cor}\label{HilbertCor}
Given any function $h\colon A\to\N$, there exists a finite set
$D\tm A$ such that
\begin{enumerate}
\item every monomial ideal with Hilbert function $h$ is generated by
monomials of degree belonging to $D$, and
\item if $I$ is a monomial ideal such that $h_I(a)=h(a)$ for all $a\in D$,
then $h_I(a)=h(a)$ for all $a\in A$. \qed
\end{enumerate}
\end{cor}

\subsection*{Noetherianity of the set of Hilbert functions}
Let $X=\{X_1,X_2,\dots\}$ be a countably infinite set of indeterminates.
For $n\geq 1$ we write
$X^{\langle n\rangle}=\{X_1,\dots,X_n\}^\diamond$, a subset of $X^\diamond$.
We identify the set $X^{\langle n\rangle}$, ordered by divisibility $\leq$, 
with 
the set $\N^n$, ordered by the product ordering, in the usual way
(see Section~\ref{Orderings-Section}).
In our last application, we will be concerned with the grading
of $X^{\langle n\rangle}$ by $A=\N$ given by
$\delta(X^\nu)=\nu_1+\cdots+\nu_n$ for $\nu=(\nu_1,\dots,\nu_n)\in\N^n$.
(This corresponds to the usual grading of monomials in a 
polynomial ring.) 
For $n\geq 1$
let $\cal H_n$ be the set of all functions $h\colon\N\to\N$ which arise
as a Hilbert function $h=h_F$ for  some
$F\in{\cal F}(X^{\langle n\rangle},\leq)$, and put $\cal H:=\bigcup_{n\geq 1}
\cal H_n$. We consider $\cal H$
as an ordered set via the product ordering: $h\leq h'$ if and only if
$h(a)\leq h'(a)$ for all $a\in\N$. We then have:

\begin{prop}\label{H-is-Noetherian}
$\cal H$ is Nash-Williams.
\end{prop}

This proposition generalizes Corollary~3.4 in \cite{AschenbrennerPong}:

\begin{cor}
$\cal H_n$ is Noetherian, for every $n\geq 1$. \qed
\end{cor}

The applications above (e.g., Proposition~\ref{KeyProp}) 
just needed Maclagan's finiteness principle,
i.e., $\N^n$ is strong\-ly Noetherian.
Proposition~\ref{H-is-Noetherian} 
is a consequence of a deeper result of Nash-Williams' theory, of
which we only state the special case needed here. 

\begin{theorem}\rom{(Nash-Williams, \cite{NW3})}\label{Deep}
If $(S,\leq)$ is a Nash-Williams ordered set, then
$(S^\diamond,\leq^\diamond)$ is Nash-Williams.
\end{theorem}

We totally order $X$ by $\preceq$ so that $X_1 \prec X_2 \prec\cdots$. 
By the previous theorem,
$\bigl({\cal F}(X^\diamond,\preceq^\diamond),\supseteq\bigr)$ is Noetherian.
We denote 
the restriction of the ordering $\preceq^\diamond$ on $X^\diamond$
to  $X^{\langle n\rangle}$ also by $\preceq^\diamond$. 
We have
$${\cal F}\bigl(X^{\langle n\rangle},\preceq^\diamond\bigr) \subseteq
{\cal F}\bigl(X^{\langle n\rangle},\leq\bigr).$$
The monomial ideals $F\in{\cal F}\bigl(X^{\langle n\rangle},\preceq^\diamond\bigr)$ 
are commonly called
{\bf strongly stable.}
It is a well-known consequence of a theorem of Galligo
(see, e.g., \cite{Eisenbud}, Section~15.9) that given a
final segment $F$ of $\bigl(X^{\langle n\rangle},\leq\bigr)$
there exists a strongly stable final segment of
 $\bigl(X^{\langle n\rangle},\leq\bigr)$ with the same Hilbert
function $h_F$.
If $F\in{\cal F}\bigl(X^{\langle n\rangle},\leq\bigr)$
and $F'$ denotes the final segment of $(X^{\langle n+1\rangle},\leq)$
generated by $F$, then $h_F \leq h_{F'}$, since $h_{F'}(m)=\sum_{i=0}^m h_F(i)$
for all $m$. We can now prove Proposition~\ref{H-is-Noetherian}:

\begin{proof}[Proof \textup{(Proposition~\ref{H-is-Noetherian})}]
Let $f\colon B\to {\cal H}$ be an $\cal H$-array, say $f(b)=h_{F_b}$ with
$F_b\in {\cal F}\bigl(X^{\langle n_b\rangle}, \leq\bigr)$, $n_b\geq 1$ 
for all $b$.  We need to find
$b\tri c$ such that $h_{F_{b}}\leq h_{F_{c}}$.
After passing to a subarray of $f$
if necessary, we may assume that $n_{b} \leq
n_{c}$ for all $b\tri c$ in $B$. (By Lemma~\ref{BadOrPerfect}.)
By the remarks above, we may further assume
that each $F_b$ is strongly stable, that is, $F_b\in
{\cal F}\bigl(X^{\langle n_b\rangle},\preceq^\diamond\bigr)$ for all $b$.
By Theorem~\ref{Deep}, there exist $b\tri c$ such that 
$(F_{b}) \supseteq (F_{c})$
in ${\cal F}\bigl(X^{\diamond},\preceq^\diamond\bigr)$, where $(F_{b})$ 
and $(F_{c})$
denote the final segments of $(X^{\diamond},\preceq^\diamond)$ generated
by $F_{b}$ and $F_{c}$, respectively. Now 
$F_{b}':=(F_{b})\cap X^{\langle n_{c}\rangle}$
is the final segment of $(X^{\langle n_{c}\rangle},\preceq^\diamond)$
generated by $F_{b}$. We have $F_{b}'\supseteq F_{c}$,
in ${\cal F}\bigl(X^{\langle n_{c}\rangle},\preceq^\diamond\bigr)$,
hence $h_{F_{b}} \leq h_{F_{b}'}\leq h_{F_{c}}$ as required.
\end{proof}

\part{Multi-stage Stochastic Integer Pro\-gram\-ming}

\section{Preliminaries: Test Sets}\label{Preliminaries: Test Sets}

\noindent
For a given matrix $A\in\Z^{l\times d}$, where $d,l\in\N$, $d,l>0$, consider the family of op\-ti\-mi\-za\-tion
problems
\begin{equation}\label{Family of problems} \tag{$\operatorname{IP}_{b,c}$}
\min\bigl\{c^\intercal z: Az=b,\ z\in\N^d\bigr\}
\end{equation}
as $b\in\R^l$ and $c\in\R^d$   vary. We call $A$ the {\bf coefficient matrix}
of this family (or of a particular instance of it).
One way to solve such a problem for given $c$ and $b$
is to start with a feasible solution
$z$ (i.e., a vector $z\in\N^d$ such that $Az=b$)
and to replace it by another feasible solution $z-v$ with
smaller objective value $c^\intercal (z-v)\in\R$,
as long as we find such a vector $v\in\Z^d$ that
improves the current feasible solution $z$.
Such a vector $v$ is called an {\bf improving vector} for $z$.
If the problem instance \eqref{Family of problems}
is solvable, this augmentation process has
to stop (with an optimal solution). Note that for given $b$ and $c$ and
any feasible solution $z$ of \eqref{Family of problems},
a vector $v\in\Z^d$ is an improving vector for $z$ if and only if the
following three conditions are satisfied:
\begin{enumerate}
\item $Av=0$,
\item $v\leq z$ (in the product ordering on $\Z^d$), and
\item $c^\intercal v>0$.
\end{enumerate}
The key step in this scheme is to find such an  improving
vector $v$. Test sets provide such
vectors: a subset of $\Z^d$
is called a {\bf universal test set} for the
above problem family if it contains, for any given choice of $c\in\R^d$
and $b\in\R^l$, an improving vector for any non-optimal feasible
solution $z$ to the given specific problem.
Clearly, if we have a finite universal test set, we can easily find an
improving vector in the augmentation procedure.

The notion of a universal test set was introduced by Graver  \cite{Graver:75}
in $1975$. He also gave a simple construction of a finite
universal test set.
The {\bf Graver basis $G(A)$} associated to $A$ consists
of all $\red$-minimal {\it non-zero}\/ integer solutions to $Az=0$.
Here $\red$ is the ordering of $\Z^d$ defined in
Section~\ref{Orderings-Section}. (Note that $0$ is the only
$\red$-minimal solution to $Az=0$.) As we have
seen above, $G(A)$ is always finite. (Corollary \ref{red-Noetherian}.)
Moreover, $G(A)$ is symmetric: if $v\in G(A)$ then also $-v\in G(A)$.
Graver showed  that $G(A)$ is indeed a universal
test set for the above problem family. This leads to the following algorithm
to compute,  uniform in $b$ and $c$,
an optimal solution
to  \eqref{Family of problems} from a feasible one:

\begin{algorithm}\label{Augmentation algorithm}
(Augmentation algorithm)

\smallskip

\Kw{Input:} a feasible solution $z$ to \eqref{Family of problems} and
a finite set $G\subseteq \Z^d$ containing $G(A)$.

\smallskip

\Kw{Output:} an optimal solution to \eqref{Family of problems}.

\smallskip

\Kw{while} there is some $v\in G$ such that $c^\intercal v>0$ and
$v\leq z$  \Kw{do}

\hskip 2.5em $z:=z-v$

\Kw{return} $z$
\end{algorithm}

\noindent
How does one find {\em some}\/ feasible solution to \eqref{Family of problems}
to begin with? Universal test sets can
also be employed to find an initial feasible solution by a
construction
similar to Phase I in the simplex method for linear
optimization.

\begin{notations}
For
$a\in\Z$ we put $$a^+:=\max\{a,0\}\in\N,\qquad a^-:=\max\{-a,0\}\in\N,$$ and
for
$z=(z_1,\dots,z_d)\in\Z^{d}$, we put
$$z^+:=(z_1^+,\dots,z_d^+)\in\N^d,\qquad
z^-:=(z_1^-,\dots,z_d^-)\in\N^d,$$
so $z=z^+-z^-$. We also let $c(z)$ denote the vector
in $\Z^d$ whose $i$-th component is $-1$ if $z_i > 0$ and $0$ otherwise.
\end{notations}

Slightly modifying Algorithm~\ref{Augmentation algorithm}
yields the following algorithm (for whose  termination and correctness see
 \cite{Hemmecke:PSP,Hemmecke:Diss}):

\begin{algorithm}\label{Finding a feasible solution}
(Finding a feasible solution)

\smallskip

\Kw{Input:} a solution $z\in\Z^d$ to \eqref{Family of problems} and
a finite set $G\subseteq \Z^d$ containing $G(A)$.

\smallskip

\Kw{Output:} a feasible solution to \eqref{Family of problems}, or
``FAIL'' if no such solution exists.

\smallskip

\Kw{while} there is some $v\in G$ such that $c(z)^\intercal v>0$ and
$v\leq z^+$  \Kw{do}

\hskip 2.5em $z:=z-v$

\Kw{if} $z\geq 0$ \Kw{then} \Kw{return} $z$ \Kw{else} \Kw{return} ``FAIL''
\end{algorithm}

\noindent
Graver proved finiteness of $G(A)$; however, he did not give an
algorithm to compute $G(A)$ from $A$. The following simple
completion procedure, due to Pottier \cite{Pottier:96},
solves this problem. We write $\ker(A):=\{z\in\Z^d:Az=0\}$ (a $\Z$-submodule
of $\Z^d$).

\begin{algorithm} \label{Completion Algorithm}
(Completion procedure)

\smallskip

\Kw{Input:} a finite symmetric generating set for
the $\Z$-module $\ker(A)$.

\smallskip

\Kw{Output:} a finite subset $G$ of $\Z^d$ containing $G(A)$.

\smallskip

$G:=F \setminus \{0\}$

$C:=\{f+g:f,g\in G\}$

\Kw{while} $C\neq \emptyset $ \Kw{do}

\hskip 2.5em $s:=$ an element in $C$

\hskip 2.5em $C:=C\setminus\{s\}$

\hskip 2.5em $f:=\normalForm(s,G)$

\hskip 2.5em \Kw{if} $f\neq 0$ \Kw{then}

\hskip 5em $C:=C\cup\{f+g:g\in G\}$

\hskip 5em $G:=G\cup {\{f\}}$

\Kw{return} $G$.
\end{algorithm}

\noindent
Behind the function $\normalForm(s,G)$ there is the following algorithm,
which upon input of a finite set $G=\{g_1,\dots,g_n\}$ of
non-zero vectors in $\Z^d$, returns a vector $f$ with the property that
$s$ has a representation  as a finite sum
\begin{equation}\label{Representation of s}
s=f+\sum_i a_ig_i \qquad\text{with
$a_i\in\N$, $g_i\in G$, $a_ig_i\red s$ and $g_i\not\red f$
for all $i$.}
\end{equation}
The vector $f$ is called a {\bf normal
form} of $s$ with respect to the set $G$.
The algorithm proceeds by successively reducing
$s$ by elements of $G$:

\begin{algorithm}
(Normal form algorithm)

\smallskip

\Kw{Input:} a vector $s\in\Z^d$, a finite set $G$ of non-zero vectors in $\Z^d$.

\smallskip

\Kw{Output:} a normal form $f$ of $s$ with respect to $G$.

\smallskip

\Kw{while} there is some $g\in G$ such that $g\red s$ \Kw{do}

\hskip 2.5em $s:=s-g$

\Kw{return} $f:=s$
\end{algorithm}

\noindent
As $\|s-g\|_1<\|s\|_1$ for $0\neq g\red s$ in $\Z^d$, the algorithm always
terminates. Here $\|s\|_1=|s_1|+\cdots+|s_d|$ for
$s=(s_1,\dots,s_d)\in\Z^d$.

In the subsequent sections we will generalize the main ingredients
of this completion procedure. The objects in our algorithm, however, will
turn out to be more complicated structures than mere vectors.
We finish this section with a useful fact needed in the next section.

%
%
%


\begin{lemma}\label{AB-Lemma}
Let $d,e,l>0$ be integers,
$A\in\Z^{l\times d}$, $B\in\Z^{l\times e}$, $y\in\Z^d$, $z\in\Z^e$,
and $b:=-Bz$. If $(y,z)\in G(A|B)$,
then either $z=0$ and $y\in G(A)$, or
$z\neq 0$ and $(y,1) \in G(A|{-b})$.
\end{lemma}
\begin{proof}
Suppose that $(y,z)\in G(A|B)$; then $(y,1)\in\ker(A|{-b})$.
If $z=0$, then clearly $y\in G(A)$, so we may
assume $z\neq 0$. Let $u\in\Z^d$, $v\in\Z$
with $Au=vb$ and $(u,v)\red (y,1)$. Then either $v=0$ or $v=1$.
If
$v=1$, then $Au=b=-Bz$, so $(u,z)\in\ker(A|B)$ and $0\neq (u,z)\red (y,z)$,
which yields $u=y$ by $(y,z)\in G(A|B)$, hence $(u,v)=(y,1)$. If $v=0$ then $Au=0$, so $(y-u,z)\in\ker(A|B)$ and
$0\neq (y-u,z)\red (y,z)$, which yields $u=0$ by  $(y,z)\in G(A|B)$, hence
$(u,v)=(0,0)$. This shows that $(y,1)\in G(A|{-b})$.
\end{proof}

%

\section{Building Blocks}\label{Building Blocks}

\noindent
In this section we study multistage stochastic integer
programs. We begin by describing the basic setup, and then use
results from the previous sections (in particular
Corollary~\ref{Power-Set-Cor}) to show that
the Graver basis elements of the coefficient matrices of
stochastic integer programs of this type can be
constructed from only finitely many ``building blocks,'' as
the number $N\geq 1$ of scenarios varies.


\subsection*{Basic setup}
In the following we fix integers $k\geq 0$ and $l\geq 1$.
Let $T_0,T_1,\dots,T_k$ be a sequence of integer matrices
with each $T_s$ having $l$ rows and $n_s$ columns, where $n_s\in\N$, $n_s\geq 1$
for $s=0,\dots,k$. We set
$A_{0,N}:=(T_0)$, and we recursively define for
$s=1,\ldots,k$
\[
A_{s,N}:=\left(
\begin{array}{ccccc}
T_{s,N}      & A_{s-1,N} & 0       & \cdots & 0\\
T_{s,N}      & 0       & A_{s-1,N} & \cdots & 0\\
\vdots   & \vdots  & \vdots  & \ddots & \vdots\\
T_{s,N}      & 0       & 0       & \cdots & A_{s-1,N}\\
\end{array}
\right)
\]
where $T_{s,N}=\left(\begin{smallmatrix} T_s \\ \vdots \\ T_s\end{smallmatrix}\right)$
consists of $N^{s-1}$ many copies of the matrix $T_s$.

\begin{example}
For $k=2$ and $N=2$ we have $A_{0,2}=(T_0)$, $A_{1,2}=
\left(\begin{smallmatrix}
T_1 & T_0  &     \\
T_1 &      & T_0
\end{smallmatrix}\right)$ and
$A_{2,2}=
\left(\begin{smallmatrix}
T_2 & T_1 & T_0 &     &     &     &     \\
T_2 & T_1 &     & T_0 &     &     &     \\
T_2 &     &     &     & T_1 & T_0 &     \\
T_2 &     &     &     & T_1 &     & T_0
\end{smallmatrix}\right)$.
\end{example}

As discussed in the introduction, matrices having the form of
$A_{k,N}$ arise as the coefficient matrices
of $(k+1)$-stage stochastic integer programs.

\begin{remarks}
Note that we assume here that the scenario tree of our stochastic
optimization problem
splits into exactly $N$ subtrees at every stage.
This simplifying condition can
easily be achieved by introducing additional scenarios with
vanishing conditional probabilities. Also note that
the coefficient matrix $A_{1,N}$ for a two-stage stochastic
integer program with $N$ scenarios differs somewhat from the
general form of a coefficient matrix given in the introduction:
in the description of a two-stage stochastic program
\[
\min\left\{c^\intercal x+\sum_{i=1}^N \pi_i\cdot(d^\intercal
y_i):Ax=a,\ Wy_i=h_i-Tx,\ x\in\N^m,\ y_i\in\N^n\right\}
\]
as in the introduction,
we read the constraints $Ax=a$ as $Ax+0y_i=a$, and in doing so, we may
rewrite the problem matrix as
\[
\left(
\begin{array}{cccc}
T' & W' & & \\
\vdots & & \ddots & \\
T' &   & & W' \\
\end{array}
\right)
\]
with $T'=\left(\begin{smallmatrix} A\\T\end{smallmatrix}\right)$ and
$W'=\left(\begin{smallmatrix} 0\\W\end{smallmatrix}\right)$. Thus,
we can safely avoid stating the constraints $Ax=a$ explicitly. The
same holds true in the multi-stage situation, as we can apply a
similar reformulation of the problem.
\end{remarks}

For $s=0,\dots,k$, the matrix $A_{s,N}$ has
$N^s\cdot l$ rows and
\[
d_{s,N}=n_s+Nd_{s-1,N}=n_s+Nn_{s-1}+N^2d_{s-2,N}=\cdots=\sum_{i=0}^s
N^{s-i}n_i
\]
columns. In particular, the full $(k+1)$-stage problem matrix $A_{k,N}$
has
\[
d_{k,N}=\sum_{i=0}^k N^{k-i}n_i
\]
columns, which corresponds to the number of variables of the
corresponding stochastic integer problem.
In the following $h$ and $s$ range over the set
$\{0,\dots,k\}$. We also put $$n(h) := n_0+\cdots+n_h\quad (=d_{h,1}).$$

\subsection*{Building blocks}
Let $N\geq 1$ and $[1,N]:=\{1,2,\ldots,N\}$. For $s=0,\dots,k$
we fix an enumeration
of the set $[1,N]^s$ of sequences of elements from $[1,N]$ of length $s$,
say in lexicographic order. By convention
$[1,N]^0=\{\varepsilon\}$ where $\varepsilon=(\ \ )$ denotes the empty sequence.
We write $$z=\bigl(z_{\alpha}:
\alpha\in [1,N]^s,\ s=0,\dots,h\bigr)$$ to denote a vector $z$ consisting of the
$1+N+N^2+\cdots+N^h$ many integer vector entries $z_\alpha\in\Z^{n_{h-s}}$, where
$\alpha\in [1,N]^s$ are indexed according to the enumeration of $[1,N]^s$.
For $s=0,\dots,h$ and $\alpha\in [1,N]^{s}$, we call $z_{\alpha}$
a {\bf building block} of height $s$ of $z$.
We have
$z\in\ker(A_{h,N})$ if and only if for all $\alpha=(\alpha_1,\dots,\alpha_h)\in [1,N]^{h}$
\[
T_hz_{\alpha|0}+T_{h-1}z_{\alpha|1}+\cdots+
T_0z_{\alpha|h}=\sum_{i=0}^h T_iz_{\alpha|h-i}=0,
\]
where $\alpha|i:=(\alpha_1,\dots,\alpha_i)$ for $i=0,\dots,h$.
(Here again $\Z^0=\{\varepsilon\}$, so
$\alpha|0=\varepsilon$ for all $\alpha\in\Z^h$.)
This relation can also be
written as
\[
A_{h,1}
\left(
\begin{array}{c}
z_{\alpha|0}\\
z_{\alpha|1}\\
\vdots\\
z_{\alpha|h}
\end{array}
\right)=
0.
\]

\subsection*{Vector-trees}
\begin{figure}[tbh]
\begin{center}
\epsfig{file=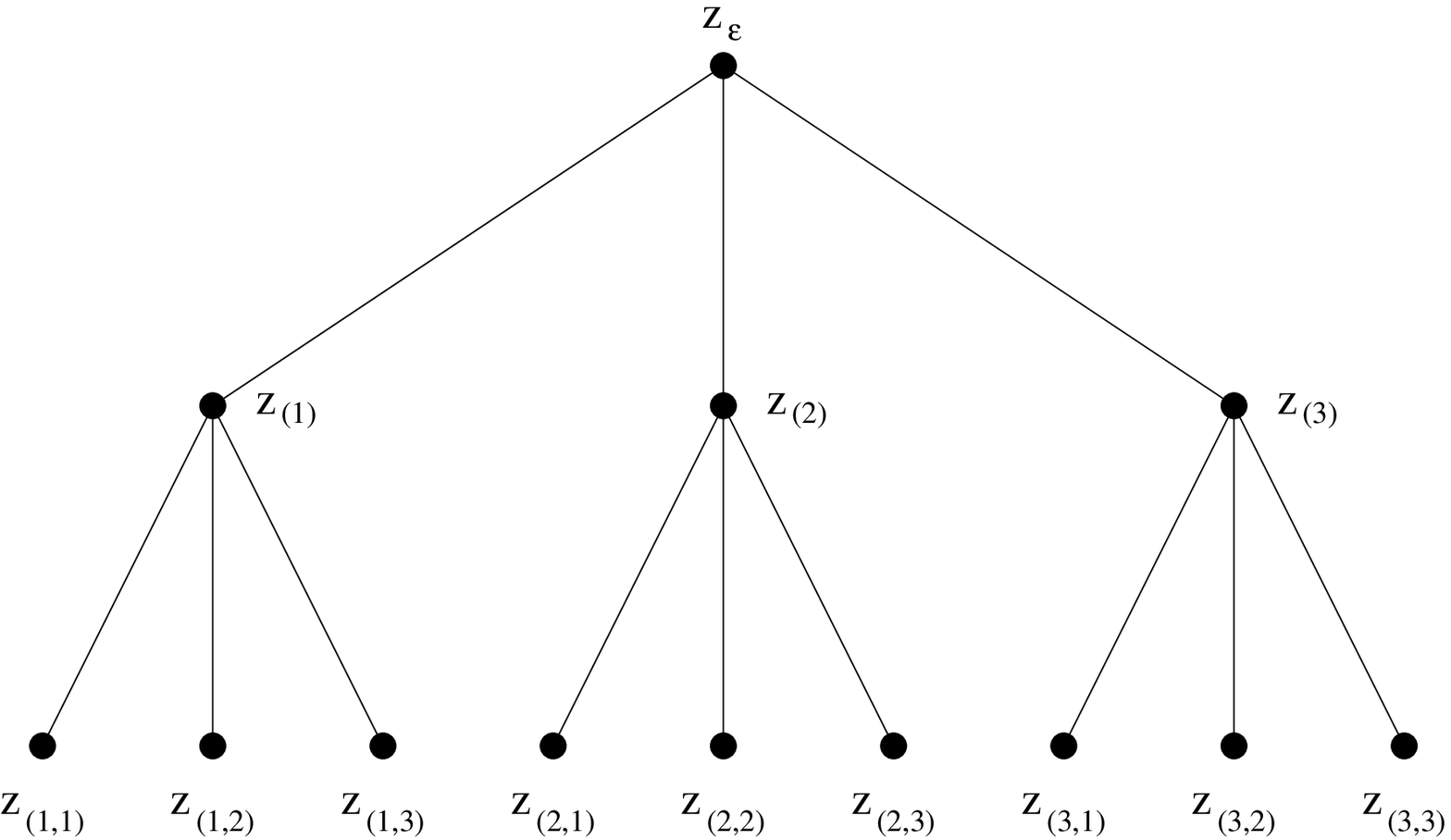, height=7cm}\\[.3cm]
\refstepcounter{figure} \label{Vector-Tree} Figure \thefigure: Representing
multi-stage solutions as trees
\end{center}
\vspace*{-0.3cm}
\end{figure}

It is useful to visualize a vector $z\in\Z^{d_{h,N}}$ as a
(directed, labeled, unordered) tree of height $h$
whose nodes are labeled by its building blocks
$z_{\alpha}$, and there is an edge from the node labeled by $z_{\alpha}$ to
the node labeled by
$z_{\beta}$ exactly if  $\alpha$ is an initial segment of
$\beta$, that is, $\alpha=\beta|s$ for some
$s=0,\dots,h$. (See Figure~\ref{Vector-Tree}.) We call this tree
the {\bf vector-tree} $\cal T(z)$ associated to $z$.
Such a vector-tree associated to  an
element of $\Z^{d_{h,N}}$ is a full $N$-ary tree (i.e.,
every internal node has exactly $N$ children).
In the following we will consider
trees of a similar structure which may split into possibly infinitely
many sub-trees at each stage:

\begin{definition}
A {\bf vector-tree} of height $h$
is a (non-empty, countable, directed, labeled)
tree which is balanced of height $h$
(i.e., every path from the root to a leaf has the same length $h$) and
whose nodes of height $s$ are labeled by integer vectors
in $\Z^{n_{h-s}}$, for $s=0,\dots,h$.
(Sometimes we also call the labels of a vector-tree its
{\bf building blocks}.)
We denote the label of the root
of a vector-tree $\cal T$ by $\root({\cal T})\in\Z^{n_{h}}$.
\end{definition}

We have an obvious notion of
isomorphism of vector-trees (as an isomorphism of directed graphs
which preserves the labeling). In the following we identify isomorphic
vector-trees; this permits us to speak of ``the set of all vector-trees''.
For every full $N$-ary vector-tree $\cal T$
of height $h$ there exists an
element $z$ of $\Z^{d_{h,N}}$ whose associated vector-tree is $\cal T$;
we have $\root\bigl(\cal T(z)\bigr)=z_{\varepsilon}$.

\begin{definition}
If ${\cal S}$ and ${\cal T}$ are vector-trees (of possibly different
heights), we
say that $\cal S$ is a {\bf sub-vector-tree} of $\cal T$ if
\begin{enumerate}
\item ${\cal S}$ is a labeled subtree of
${\cal T}$, that is, the underlying graph of $\cal S$ is
a subgraph of the underlying graph of $\cal T$,
and the labeling of
the nodes of $\cal S$ agrees with their corresponding labeling in $\cal T$;
\item ${\cal S}$ is closed downwards in $\cal T$, that is, if there is a
path in $\cal T$ from the root of $\cal S$ to a node $a$ in $\cal T$, then
$a$ is a node of $\cal S$.
\end{enumerate}
If ${\cal S}$ is a sub-vector-tree of $\cal T$ where $\cal S$ has
$h-1$ and $\cal T$ has height $h$, then $\cal S$ is called an
{\bf immediate} sub-vector-tree of $\cal T$.
\end{definition}

\begin{example}
In Figure~\ref{Vector-Tree}, the labeled subtree consisting of the nodes
labeled by $z_{(1)}$, $z_{(1,1)}$, $z_{(1,2)}$, $z_{(1,3)}$
is a sub-vector-tree
of the vector-tree associated to $z$, whereas the labeled subtree consisting
of the nodes labeled by $z_{\varepsilon}$, $z_{(1)}$, $z_{(1,1)}$,
$z_{(1,2)}$ is not.
\end{example}

\subsection*{Paths in vector-trees}
Every path in a vector-tree $\mathcal T$
of height $h$ from its root to one of it leaves
is called a {\bf maximal path.} Every maximal path in a vector tree of
height $h$ has length $h$. If $v_0,v_1,\dots,v_h$ are the successive labels
of the nodes on a path $P$ in $\cal T$
we also say that $P$ is {\bf labeled by} $(v_0,\dots,v_h)$, and we
call $(v_0,\dots,v_h)\in\Z^{n(h)}$ the {\bf label of $P$.}
We write $\paths(\cal T)$ for the set of labels of maximal paths in
$\cal T$.

\begin{definition}
We say that a vector-tree $\cal T$ of height $h$ has {\bf value} $b\in\Z^l$ if
for every  $(v_0,v_{1},\ldots,v_h)\in\paths(\cal T)$ we have
\[
T_h v_0+T_{h-1}v_{1}+\cdots+ T_0v_h=b.
\]
\end{definition}

Note that if ${\cal T}$ is a vector-tree of height $h > 0$ and value $b$,
then each immediate sub-vector-tree of $\cal T$
has value $b-T_h \root({\cal T})$.
An element $z$ of $\Z^{d_{h,N}}$ lies in $\ker(A_{h,N})$ if and only if
$\cal T(z)$ has value $0$.

\begin{definition}
Let $\cal T$ be a vector-tree of height $h$ and $N$ a positive integer.
A vector $z\in\Z^{d_{h,N}}$ with the property
that $\paths\big(\cal T(z)\big)\subseteq\paths(\cal T)$ is said to be
{\bf constructible from $\cal T$.} We denote by
$\langle {\cal T}\rangle_{N}$ the set of vectors in $\Z^{d_{h,N}}$ which
are constructible from $\cal T$.
\end{definition}

\begin{remarks}
We have $\langle\cal T\rangle_1 = \paths(\cal T)$.
If $\cal T$ has value $b$, then for every $z\in
\langle {\cal T}\rangle_{N}$, the vector-tree $\cal T(z)$
has value $b$.
\end{remarks}

We say that a vector-tree $\cal T$ is {\bf tight} if
the labels of the children of internal nodes of $\cal T$ are
pairwise distinct. For example, the vector-tree
$$
\addtolength{\varrowlength}{-20pt}
\addtolength{\harrowlength}{-15pt}
\addtolength{\sarrowlength}{-45pt}
\commdiag{  & 0 & \cr
            &\arrow(-1,-1)\arrow(1,-1) \cr
          0 &   & 2
          }
$$
is tight, whereas
$$
\addtolength{\varrowlength}{-20pt}
\addtolength{\harrowlength}{-15pt}
\addtolength{\sarrowlength}{-45pt}
\commdiag{   & 0 & \cr
            &\arrow(-1,-1)\mapdown \arrow(1,-1) \cr
          0 &  0 & 2
          }
$$
is not. The following is easy to show:

\begin{lemma}
Let $S$ be a non-empty subset of $\Z^{n(h)}$ with $z_{\epsilon}=r$
for all $z\in S$. There exists a unique tight vector-tree ${\cal T}(S)$
of height $h$ \textup{(}up to isomorphism\textup{)}
with  $S=\paths\big(\cal T(S)\big)$.
\end{lemma}

In particular, for every vector-tree $\cal S$ of height $h$
there exists a unique tight vector-tree $\cal T$ of the same height with
$\paths(\cal T)=\paths(\cal S)$, namely $\cal T=\cal T\big(\paths(\cal S)\big)$.
Note also that in the context of the previous lemma, if
$S\subseteq \ker(A_{h,1})$, then $\cal T(S)$ has value $0$.

\subsection*{Reducibility of vector-trees}
Next,
we define a ``reducibility relation'' between vector-trees of the same height.
This relation bears a formal resemblance to Milner's
``simulation quasi-ordering'' for transition systems  \cite{Milner-R}.

\begin{definition}
Let ${\cal S}$ and ${\cal T}$ be vector-trees of height  $h$.
\begin{enumerate}
\item For $h=0$ we let
${\cal S}\red_0 {\cal T}$ if $\root({\cal S}) \red
\root({\cal T})$ (in $\Z^{n_0}$).
\item For $h>0$ we let ${\cal S}\red_h {\cal T}$ if
$\root({\cal S})\red \root({\cal T})$  (in $\Z^{n_h}$) and
for every immediate sub-vector-tree ${\cal T}'$ of ${\cal T}$
there is an immediate
sub-vector-tree ${\cal S}'$ of ${\cal S}$ such that
${\cal S}'\red_{h-1} {\cal T}'$.
\end{enumerate}
\end{definition}

There
is an obvious algorithm to decide, given {\it finite}\/
vector-trees $\cal S$ and $\cal T$ of height $h$, whether $\cal S\red_h\cal T$.
Clearly $\red_h$ is a
quasi-ordering  on the set of vector-trees of
height $h$.
%
%
In the proof of the next lemma
we relate $\red_h$ with the quasi-ordering $\red_h$
on $\cal P_h(\Z^{n(h)})$ defined
after Theorem~\ref{NmBetterNoeth};
we freely use the notations introduced there.

\begin{lemma}\label{red_h is Noetherian}
The quasi-ordering $\red_h$ on the set of vector-trees of height $h$
is Noetherian.
\end{lemma}
\begin{proof}
First we put $${\cal P}_0^*(\Z^n) := \Z^n, \quad
{\cal P}_{h}^*(\Z^n) := {\cal P}\bigl({\cal P}_{h-1}^*(\Z^n)\bigr)\setminus\{\emptyset\}\text{ for $h>0$.}$$
So ${\cal P}_h^*(\Z^n)$ is a subset of ${\cal P}_h(\Z^n)$, for each $h$
and $n>0$.
By induction on $h= 0,\dots,k$ we now define for all $m,n\in\N$, $m,n>0$,
an operation $$(v,X)\mapsto v\ast_h X\colon
\Z^{m}\times {\cal P}_{h}^*(\Z^{n}) \to
{\cal P}_{h}^*(\Z^{m+n})$$
as follows. For $h=0$ we have
${\cal P}_{h}^*(\Z^{n})=\Z^n$;
we let $v\ast_0 w\in\Z^{m+n}$ be the concatenation of
$v\in\Z^m$ and $w\in\Z^n$.
For $h>0$ we let $$v\ast_k X := \bigl\{ v\ast_{h-1} X' : X'\in X\bigr\}
\quad\text{for
all $v\in\Z^{m}$, $X\in{\cal P}_{h}^*(\Z^{n})$.}$$
Note that for $v,w\in\Z^{m}$ and
$X,Y\in {\cal P}_{h}^*(\Z^{n})$ we have: if
$v\ast_h X\red_{h} w\ast_h Y$ in ${\cal P}_{h}\bigl(\Z^{m+n}\bigr)$,
then $v\red w$ and
$X\red_{h} Y$ (in ${\cal P}_{h}(\Z^{n})$).

Next, we define a map $\varphi_h$ which associates to every
vector-tree $\cal S$ of height $h$ an element $\varphi_h({\cal S})$
of ${\cal P}_{h}^*(\Z^{n(h)})$
as follows: For $h=0$ put $\varphi_0({\cal S}):=\root({\cal S})\in
\Z^{n_0}$. For $h>0$
let $$\varphi_h({\cal S}) := \bigl\{ \root({\cal S})\ast_h
\varphi_{h-1}({\cal S}') :
\text{${\cal S}'$ sub-vector-tree of $\cal S$ of height $h-1$}\bigr\}.$$
By induction on $h$
it is easy to verify that $\varphi_h$ is a quasi-embedding
of the set of vector-trees of height $h$, quasi-ordered by $\red_h$ as defined
above, into
${\cal P}_{h}^*\bigl(\Z^{n(h)}\bigr)$, quasi-ordered by the restriction of the
quasi-ordering $\red_h$ of ${\cal P}_{h}\bigl(\Z^{n(h)}\bigr)$
defined after Theorem~\ref{NmBetterNoeth}. The lemma now follows from
Corollary~\ref{Power-Set-Cor}.
\end{proof}

\begin{cor} \label{Sequence of non-reducible trees is finite}
There is no infinite sequence $\cal T_1,\cal T_2,\ldots$ of
vector-trees of height $h$ with $\cal T_i\not\red_{h}
\cal T_j$ whenever $i<j$. \qed
\end{cor}


Given vector-trees $\cal S$ and $\cal T$ of height $h$ and an integer $N\geq 1$
we write ${\cal S} \red_{h,N} {\cal T}$ if
for every $z\in\langle {\cal T}\rangle_{N}$ there exists
$y\in\langle {\cal S}\rangle_{N}$ such that $y\red z$
(in $\Z^{d_{h,N}}$). We put ${\cal S}\red_{h,\infty} {\cal T}$ if
${\cal S} \red_{h,N} {\cal T}$ for all $N\geq 1$.

\begin{lemma}\label{Tree reducibility implies vector reducibility}
If ${\cal S}\red_h {\cal T}$, then ${\cal S} \red_{h,\infty} {\cal T}$.
\end{lemma}

\boproof
Suppose that  ${\cal S}\red_h {\cal T}$, and
let $N\geq 1$. By induction on $h=0,\dots,k$ we show that
given $z\in\langle {\cal T}\rangle_{N}$ we can construct a
vector $y\in\langle {\cal S}\rangle_{N}$ such that
$y\red z$. Suppose first that $h=0$. Then
$z=\root({\cal T})$, and $\root({\cal S})\red \root({\cal T})$.
Hence for $y$ we may take $y=\root({\cal S})$.
Now assume that $h>0$, and let $z_1,\dots,z_N\in\Z^{d_{h-1,N}}$ such that
$\cal T(z_1),\dots,\cal T(z_N)$ are the immediate sub-vector-trees
of $\cal T(z)$. Then for each $i=1,\dots,N$ there is
an immediate sub-vector-tree ${\cal T}_i$ of $\cal T$  with $z_i\in
\langle{\cal T}_i\rangle_{N}$. Since ${\cal S}\red_h{\cal T}$
we have $\root({\cal S})\red \root({\cal T})$, and for every $i$ there
exists a sub-vector-tree ${\cal S}_i$ of $\cal S$ of height $h-1$ with
${\cal S}_i \red_{h-1} {\cal T}_i$. By induction hypothesis
there exist $y_1,\dots,y_N\in\langle {\cal S}_i\rangle_{N}$
such that $y_i \red z_i$ in $\Z^{d_{h-1,N}}$.
Then $y\red z$ for a suitable $y\in\Z^{d_{h,N}}$ whose vector-tree
$\cal T(y)$ has immediate
sub-vector-trees $\cal T(y_1),\dots,\cal T(y_N)$ and root
labeled by $\root({\cal S})$.
\eoproof

Lemmas~\ref{red_h is Noetherian} and
\ref{Tree reducibility implies vector reducibility} imply:

\begin{cor} \label{Tree reducibility implies vector reducibility, cor}
For every $h$, the quasi-ordering $\red_{h,\infty}$ on
the set of vector-trees of height $h$ is Noetherian. \qed
\end{cor}

\subsection*{Building blocks of Graver bases}
Recall that $G(A_{k,N})$ denotes the Graver basis of the matrix
$A_{k,N}$ (a finite subset of $\Z^{d_{k,N}}$).
%
We decompose each vector $z\in G(A_{k,N})$ into its building blocks
$z_{\alpha}$ as described above and put
$${\cal H}_{k,N}^s :=  \bigl\{ z_{\alpha} : z\in G(A_{k,N}),\
\alpha\in [1,N]^{s}\bigr\}\subseteq\Z^{n_{k-s}}\qquad\text{for $s=0,\dots,k$.}$$
We form the union
${\cal H}^s_{k,\infty}:=\bigcup_{N=1}^\infty {\cal H}_{k,N}^s$
and define
\[
{\cal H}_{k,\infty}:=\bigcup\limits_{s=0}^{k} {\cal H}_{k,\infty}^s\quad
\text{(disjoint union),}
\]
the set of building blocks of Graver bases of the matrices
$A_{k,N}$ obtained by varying $N$ over the set of positive integers.
For $k=0$ we have ${\cal H}_{0,\infty}=G(T_0)$ and thus, ${\cal
H}_{0,\infty}$ is finite. Finiteness of ${\cal H}_{1,\infty}$ was
shown in \cite{Hemmecke+Schultz:dec2SIP}. In what follows, we prove:

\begin{prop}\label{Hkinfty is finite}
The set ${\cal H}_{k,\infty}$ is finite for every $k$.
\end{prop}

Before we give the proof of this proposition,
we establish a few auxiliary facts.
We first combine the elements of ${\cal H}_{k,\infty}$ into
({\it a priori}\/ possibly infinite) vector-trees:
Given $r\in{\cal H}_{k,\infty}^0$,
we put $\cal T(r) := \cal T\big(S(r)\big)$ where
$$S(r) := \bigcup \big\{ \paths\big(\cal T(z)\big) :
\text{$z\in G(A_{k,N})$ for some $N\geq 1$, $z_{\epsilon}=r$}
\big\} \subseteq \Z^{n(k)}.$$
The vector-tree ${\cal T}(r)$ of height $k$ has root $r$ and value $0$.
For any $r\in{\cal H}_{k,\infty}^0$ and
any $N\geq 1$ the set $\langle {\cal T}(r)\rangle_{N}$
contains every vector $z\in G(A_{k,N})$ with $z_{\varepsilon}=r$.
In particular
$G(A_{k,N})\cap \langle
{\cal T}(r)\rangle_{N}\neq\emptyset$ for some $N\geq 1$.
This yields:

\begin{cor}\label{Subtrees in Hinfty are not comparable}
The vector-trees ${\cal T}(r)$, where $r$ ranges over all non-zero elements
of ${\cal H}_{k,\infty}^0$, form a $\red_{k,\infty}$-antichain.
\end{cor}

\boproof For a contradiction
suppose that ${\cal T}(r')\red_{k,\infty} {\cal T}(r)$ for some
non-zero $r'\neq r$ in ${\cal H}_{k,\infty}^0$.
So for every $N\geq 1$ and every
$z\in\langle {\cal T}(r)\rangle_{N}$ there
exists a vector $y\in\langle {\cal T}(r')\rangle_{N}$ such
that $y\red z$. Note that
$y_{\varepsilon}=\root\bigl({\cal T}(r')\bigr)\neq
\root\bigl({\cal T}(r)\bigr)=z_{\varepsilon}$, hence $y,z\neq 0$ and $y\neq z$.
Therefore, none of the vectors constructible from
${\cal T}(r)$ is an element of a Graver basis $G(A_{k,N})$, for any
$N$. This contradicts the remark preceding the corollary.
\eoproof

By
Corollaries~\ref{Tree reducibility implies vector reducibility, cor} and
\ref{Subtrees in Hinfty are not comparable}, the set
${\cal H}_{k,\infty}^0$ is finite.
We can now prove Proposition~\ref{Hkinfty is finite}:

\begin{proof}[Proof \rom{(Proposition~\ref{Hkinfty is finite})}]
We show, by induction on $k$, that ${\cal H}_{k,\infty}$ is finite, for
every choice of matrices $T_0,\dots,T_k$
as in the beginning of this section.
We already know that ${\cal H}_{0,\infty}=G(T_0)$ is finite.
Suppose $k>0$.
We have seen above that
$\cal H_{k,\infty}^0$ is finite. Hence it suffices to show that
for each $r\in\cal H_{k,\infty}^0$, the vector-tree
${\cal T}(r)$ is finite.

Suppose first that $r=0$. By induction hypothesis, ${\cal H}_{k-1,\infty}$
is finite. Hence it is enough to show that
all labels of non-root nodes of ${\cal T}(0)$
are in ${\cal H}_{k-1,\infty}$. Let $v$ be the label of a node of
height $s>0$ of ${\cal T}(0)$.
There exists an integer $N\geq 1$, a vector
$z\in G(A_{k,N})$, and $\alpha=(\alpha_1,\dots,\alpha_s)\in [1,N]^s$ such that $z_{\varepsilon}=0$ and
$z_{\alpha}=v$. Let $z'\in\Z^{d_{k-1,N}}$ such that $\cal T(z')$ is
a sub-vector-tree of
$\cal T(z)$ and $z'_{\alpha'}=v$, where
$\alpha'=(\alpha_2,\dots,\alpha_s)$
(i.e.,  $\cal T(z')$ contains our node labeled by $v$).
Then clearly
$z'\in G(A_{k-1,N})$ and thus $v\in {\cal H}_{k-1,\infty}^s$
as desired.

Now suppose $r\neq 0$. Let
$T_{k-1}'$ be the $l\times (1+n_{k-1})$-matrix
$({-T_kr}|{T_{k-1}})$ and
put $T_s':=T_s$ for $s=0,\dots,k-2$. Define
$A_{s,N}'$ in the same way as $A_{s,N}$
at the beginning of this section, with the matrices
$T_s'$ replacing $T_s$. Let ${\cal H}'_{k-1,\infty}$ be
the set of building blocks of Graver bases of $A_{k-1,N}'$, for $N\geq 1$.
By induction hypothesis, ${\cal H}'_{k-1,\infty}$ is finite, so
it is enough to show that for every node of height $s>0$ of
${\cal T}(r)$ with label $v$, we have
$v\in{\cal H}'_{k-1,\infty}$.
Let $N\geq 1$, $z\in G(A_{k,N})$, and
$\alpha=(\alpha_1,\dots,\alpha_s)\in [1,N]^s$ with $z_{\varepsilon}=r$
and $z_{\alpha}=v$.
Let $z'\in\Z^{d_{k-1,N}}$ such that $\cal T(z')$ is
a sub-vector-tree of $\cal T(z)$
with $z'_{\alpha'}=v$, where $\alpha'=(\alpha_2,\dots,\alpha_s)$.
Then $(r,z')\in G(T_{k,N}|A_{k-1,N})$.
By Lemma~\ref{AB-Lemma}, $z'$ is
a projection of an element of $G(A_{k-1,N}')$.
Therefore the building blocks
of $z'$ are in ${\cal H}'_{k-1,\infty}$. In particular
$v\in{\cal H}'_{k-1,\infty}$ as required.
\end{proof}

%
%

The proof of the proposition above suggests a procedure
for constructing $\cal H_{k,\infty}$ for $k=0,1,2,\dots$ inductively.
In the next section we describe such an algorithm.
We finish this section by a few remarks about the choice of the
reducibility relation $\red_k$.
First note the following immediate consequence of the fact that
$\paths(\cal T)=\langle\cal T\rangle_1$ for every vector-tree $\cal T$:

\begin{lemma}\label{Motivation for refined reducibility}
For vector-trees $\cal S$, $\cal T$ of height $h$ we have:
if  $\cal S \red_{h,\infty} \cal T$, then
for every $w\in\paths(\cal T)$
there exists $v\in\paths(\cal S)$ such that $v\red w$
\textup{(}in $\Z^{n(h)}$\textup{)}.
\end{lemma}

The converses of the implications in
Lemmas~\ref{Tree reducibility implies vector reducibility} and
\ref{Motivation for refined reducibility} are
false in general, as the following two simple examples show.
In particular, this (at least partly) explains why in the algorithm
for computing $\cal H_{k,\infty}$
in the next section, we cannot simply replace
$\red_k$ by the
quasi-ordering $\leq_k$ on vector-trees of height $k$ given by
\begin{equation}\label{proviso}
\cal S\leq_k \cal T \quad :\Longleftrightarrow \quad
\begin{cases}
& \parbox{20em}{for every
$w\in\paths(\cal T)$
there is $v\in\paths(\cal S)$ such that $v\red w$,}\end{cases}
\end{equation}
whose Noetherianity
is much easier to show than the Noetherianity of the quasi-ordering $\red_k$
(cf.~Corollary~\ref{Maclagan} in the case $n=n(k)$).
In both examples  $k=2$ and $n_0=n_1=n_2=2$.


\begin{example}
The vector-trees
$$
\addtolength{\varrowlength}{-20pt}
\addtolength{\harrowlength}{-15pt}
\addtolength{\sarrowlength}{-45pt}
\commdiag{
           &(0,0)                       &               & \qquad\qquad\qquad\qquad\qquad&      &(0,0)                     & \cr
           &\arrow(-1,-1) \arrow(1,-1)  &               &                               &      &\mapdown                  &\cr
(0,0)      &                            &  (0,0)        &                               &      &(0,0)                     &\cr
\mapdown   &                            & \mapdown      &                               &      &\arrow(-1,-1) \arrow(1,-1)&\cr
(0,1)      &                            & (1,0)         &                               &(0,1) &                          & (1,0) \cr
           & \cal S                     &               &                               &      & \cal T                   &
}
$$
show that we may have $\cal S\red_{2,\infty}\cal T$ and ${\cal S}\not\red_2{\cal T}$.
\end{example}

\begin{example}
Consider the vector-trees
$$
\addtolength{\varrowlength}{-20pt}
\addtolength{\harrowlength}{-15pt}
\addtolength{\sarrowlength}{-45pt}
\commdiag{
           &(0,0)                       &               & \qquad\qquad\qquad\qquad\qquad&      &(0,0)                     & \cr
           &\arrow(-1,-1) \arrow(1,-1)  &               &                               &      &\mapdown                  &\cr
(0,0)      &                            &  (0,1)        &                               &      &(0,1)                     &\cr
\mapdown   &                            & \mapdown      &                               &      &\arrow(-1,-1) \arrow(1,-1)&\cr
(1,0)      &                            & (0,1)         &                               &(1,0) &                          & (0,1) \cr
           & \cal S                     &               &                               &      & \cal T                   &
}
$$
We have $\cal S\leq_k \cal T$ (as defined in \eqref{proviso}).
Suppose that $z\in\Z^{14}$ has associated
vector-tree
$$
\addtolength{\varrowlength}{-20pt}
\addtolength{\harrowlength}{-15pt}
\addtolength{\sarrowlength}{-45pt}
\commdiag{
      &                           &               & (0,0) &               &                             &           \cr
      &                           & \arrow(-1,-1) &       & \arrow(1,-1)  &                             &            \cr
      & (0,1)                     &               &       &               &  (0,1)                      & \cr
      & \arrow(-1,-1) \arrow(1,-1)&               &       &               &  \arrow(-1,-1) \arrow(1,-1) &   \cr
(1,0) &                           & (0,1)         &       &  (1,0)        &                             & (0,1) \cr
}
$$
Then $z\in\langle\cal T\rangle_2$, but there is no $y\in\langle\cal S\rangle_2$
with $y\red z$.
\end{example}

\section{Computation of Building Blocks}\label{Computation of Building Blocks}

\noindent
Our algorithm for computing a finite set of vector containing
${\cal H}_{k,\infty}$ follows the pattern of a
completion procedure, similar to Buchberger's algorithm for computing
Gr\"obner bases of ideals in polynomial rings over fields. 
Instead of with (finite sets of) polynomials, our procedure operates
with finite vector-trees. Before we describe our algorithm, we
need  to specify some crucial ingredients for this completion procedure, 
among them the input set and a notion of normal form.
The notations and conventions adapted in the last section remain in
force until the end of the paper.

\subsection*{Adding and subtracting vector trees}
We begin by defining 
operations which allow us to construct new vector-trees from
old ones.
For this we use the following notations, for
subsets $V$ and $W$ of $\Z^{m}$, $m\geq 1$:
\begin{align*}
-V     &:= \{ -v    : v\in V\}, \\
V +  W &:= \{ v + w : v\in V,\ w\in W \}, \\
V -  W &:= \{ v - w : v\in V,\ w\in W,\ w\red v\}.
\end{align*}
Note that in general $V-W\neq V+(-W)$.
In the following 
$\cal S$ and $\cal T$ range over the set of 
vector-trees of height $k$.

\begin{definition}
We put
$$
-\cal S := \cal T\big({-\paths(\cal S)}\big), \quad
\cal S + \cal T := \cal T\big(\paths(\cal S) + \paths(\cal T)\big).
$$
If  $\cal T\red_k\cal S$, then $\paths(\cal S)-\paths(\cal T)\neq\emptyset$
by Lemmas~\ref{Tree reducibility implies vector reducibility} and \ref{Motivation for refined reducibility}, and in this case we put
$$\cal S - \cal T  := \cal T\big(\paths(\cal S) - \paths(\cal T)\big).$$
\end{definition}


\begin{remarks}\ 

\begin{enumerate}
\item We have $\paths(-\cal S)=-\paths(\cal S)$ and
$\paths(\cal S+\cal T)=\paths(\cal S)+\paths(\cal T)$. 
If $\cal T\red_k\cal S$ then 
$\paths(\cal S-\cal T)=\paths(\cal S)-\paths(\cal T)$.
\item If $\cal S$ has value $a\in\Z^l$ 
and $\cal T$ has value $b\in\Z^l$, then $-\cal S$ has value $-a$,
$\cal S+\cal T$ has value $a+b$, and if in addition $\cal T\red_{k}\cal S$,
then $\cal S-\cal T$ has value $a-b$.
\item There are obvious algorithms to compute, given
finite $\cal S$ and $\cal T$, the vector-trees $-\cal S$, $\cal S+\cal T$,
and $\cal S-\cal T$ (provided $\cal T\red_k\cal S$). 
\end{enumerate}
\end{remarks}

\subsection*{Normal forms}
We say that $\cal S^*$ is
a {\bf normal form} of $\cal S$ with respect to 
a set $G$ of vector-trees of height $k$ if
\begin{enumerate}
\item $\cal T\not\red_k\cal S^*$ for all $\cal T\in G$ with $\root(\cal T)\neq 0$,
\item there exists a sequence $\cal S_0,\dots,\cal S_n$ of vector-trees of
height $k$ such that $\cal S_0=\cal S$, $\cal S_n=\cal S^*$, and
for every $i=0,\dots,n-1$ there exists $\cal T_i\in G$ with
$\cal T_i\red_k\cal S_i$ and
$\cal S_{i+1}=\cal S_i-\cal T_i$.
\end{enumerate}
Note that for $k=0$, if we identify each vector-tree of height $0$ with the
label of its root, this notion corresponds to the notion of normal form
for elements of $\Z^{n_0}$ 
introduced in Section~\ref{Preliminaries: Test Sets}.
The following algorithm computes a normal form:

\begin{algorithm}\label{Normal form algorithm}
{(Normal form algorithm)}

\smallskip

\Kw{Input:} \parbox[t]{30em}{a finite vector-tree ${\cal S}$ of height $k$ and a finite
set $G$ of finite vector-trees of height $k$.}

\smallskip

\Kw{Output:} a normal form  $\operatorname{normalForm}(\cal S,G)$ of $\cal S$ with respect to $G$.

\smallskip

\Kw{while} there is some  ${\cal T}\in G$ such that ${\cal
T}\red_{k} {\cal S}$ and $\root(\cal T)\neq 0$ \Kw{do}

\hskip 2.5em $\cal S:=\cal S-\cal T$

\Kw{return} ${\cal S}$
\end{algorithm}

\noindent
The algorithm above terminates: if $\cal T\red_k\cal S$, then $\root(\cal S-\cal T)=\root(\cal S)-\root(\cal T)
\red\root(\cal S)$, hence if in addition 
$\root(\cal T)\neq 0$, then $||\root(\cal S-\cal T)||_1<||\root(\cal S)||_1$.
Note that if every vector-tree in $G$ has value $0$, and $\cal S$ has
value $a$, then the output $\operatorname{normalForm}(\cal S,G)$ 
also has value $a$.
Algorithm~\ref{Normal form algorithm} will be employed as a sub-program in our algorithm
for computing the set $\cal H_{k,\infty}$.
In the proof of the correctness of the latter we will use the following 
crucial lemma  and its corollary below.

\begin{lemma}\label{Useful Lemma}
Let $N\geq 1$. 
\begin{enumerate}
\item If $y\in \langle\cal S\rangle_N$ and $z\in\langle\cal T\rangle_N$,
then $y+z \in \langle\cal S + \cal T\rangle_N$.
\item If $\cal T\red_{k}\cal S$ 
and $y\in  \langle\cal S\rangle_N$,
then there exists $z\in \langle\cal T\rangle_N$ such that
$z\red y$ and 
$y-z\in \langle\cal S-\cal T\rangle_N$.
\end{enumerate}
\end{lemma}
\begin{proof}
For part (1), suppose that $y\in \langle\cal S\rangle_N$ and $z\in\langle\cal T\rangle_N$, that is, 
$\paths\big(\cal T(y)\big)\subseteq\paths(\cal S)$ and
$\paths\big(\cal T(z)\big)\subseteq\paths(\cal T)$.
Then
\begin{align*}
\paths\big(\cal T(y+z)\big) &\subseteq \paths\big(\cal T(y)\big)+\paths\big(\cal T(z)\big) \\
&\subseteq \paths(\cal S)+\paths(\cal T) 
= \paths(\cal S+\cal T),
\end{align*}
hence $y+z \in \langle\cal S + \cal T\rangle_N$ as claimed.
We show (2) by induction on $k$.
The case $k=0$ being easy, suppose that $k>0$ and 
the claim holds with $k-1$ in place of $k$. 
Since $\cal T\red_k\cal S$, we have $\root(\cal T)\red \root(\cal S)=
y_{0,\epsilon}$. Consider
$y'\in\Z^{d_{k-1,N}}$ whose vector-tree $\cal T(y')$ is an immediate
sub-vector-tree of $\cal T(y)$. Then there exists a unique
immediate sub-vector-tree $\cal S'$ of $\cal S$ with $y'\in\langle\cal S'\rangle_N$.
Since $\cal T\red_k\cal S$ we find an immediate sub-vector-tree
$\cal T'$ of $\cal T$ with $\cal T'\red_{k-1}\cal S'$.
By inductive hypothesis there exists $z'\in\langle \cal T'\rangle_N$ such that
$z'\red y'$ and $y'-z'\in\langle\cal S'-\cal T'\rangle_N$.
These remarks suffice to construct a vector $z\in\Z^{d_{k,N}}$ with the
required properties.
\end{proof}

For a set $G$ of vector-trees of height $k$ and an integer $N\geq 1$ we
put $$\langle G\rangle_N := \bigcup_{\cal T\in G} \langle\cal T\rangle_N.$$
Part (2) of the last lemma (which strengthens Lemma~\ref{Tree reducibility implies vector reducibility}) immediately implies:

\begin{cor}\label{Useful Lemma, Cor}
Let $G$ be a set of vector-trees of height $k$, and let $\cal S^*$ be a normal
form of $\cal S$ with respect to $G$. For all $N\geq 1$ and 
$s\in\langle\cal S\rangle_N$ there exist 
$f\in\langle\cal S^*\rangle_N$ and $g_1,\dots,g_n\in\langle G\rangle_N$
such that
$$s=f+\sum_i g_i, \qquad\text{ $g_i\red s$ for all $i$.}$$
\end{cor}

\subsection*{Choosing an input set}
The following lemma justifies the choice of input set for 
Algorithm~\ref{Algorithm to compute H_infty} below. 

\begin{lemma}\label{Generators}
Suppose that $k>0$, and let $F$ be a set of generators for the $\Z$-submodule
$K:=\ker(A_{k,1})$ of $\Z^{n(k)}$
which contains a set of generators for the submodule
$$K_0 := 
K\cap \big(\{0\}\times \Z^{n(k-1)} \big)$$
of $K$. Then for every $N\geq 1$, the set
$$F_N := \big\langle\big\{ \cal T(v):v\in F\big\}\big\rangle_N$$
generates $\ker(A_{k,N})$.
\end{lemma}

A generating set $F$ satisfying the hypothesis of 
the lemma can be found algorithmically by standard methods (e.g., Hermite 
normal form).
In the proof of the lemma we use the following notations. 

\begin{notations}
Let $N\geq 1$, $\alpha=(\alpha_1,\dots,\alpha_k)\in [1,N]^k$.
For $z\in\Z^{d_{k,N}}$ we
put $$v(z,\alpha):=(z_{\alpha|0},\dots,z_{\alpha|k})\in\paths\big(\cal T(z)\big),$$
and for $v=(v_0,\dots,v_k)\in\Z^{n(k)}$, we denote 
by $z(\alpha,v)$  the vector $z\in\Z^{d_{k,N}}$ which satisfies,
for $s=0,\dots,k$ and $\beta\in [1,N]^s$:
$$z_{\beta} := \begin{cases}
v_s & \text{if $\beta=\alpha|s$,}\\
0   & \text{otherwise.}
\end{cases}$$
We also let $z(v)$ be the  unique element of
$\Z^{d_{k,N}}$ whose vector-tree $\cal T=\cal T\big(z(v)\big)$
satisfies $\paths(\cal T)=\{v\}$.
Clearly, for $v\in\Z^{n(k)}$ and $z\in\Z^{d_{k,N}}$:
\begin{equation}\label{impl1}
z\in\ker(A_{k,N})\quad\Rightarrow\quad v(\alpha,z)\in\ker(A_{k,1}),
\end{equation}
 and 
$$
v\in\ker(A_{k,1})\ \Longleftrightarrow\ z(\alpha,v)\in\ker(A_{k,N})
\ \Longleftrightarrow\ z(v)\in\ker(A_{k,N}).
$$
\end{notations}

\begin{proof}[Proof \textup{(Lemma~\ref{Generators})}]
Let $z\in\ker(A_{k,N})$, where $k,N\geq 1$. Fix an arbitrary
$\alpha\in [1,N]^k$ and put $v=v(\alpha,z)$; by \eqref{impl1}
there are $a_1,\dots,a_m\in\Z$ and $v_1,\dots,v_m\in F$ 
such that $$v=a_1v_1+\cdots+a_mv_m.$$
Then $z(v_i)\in F_N$ for each $i$, and 
\begin{equation}\label{equ1}
z(v)=a_1z(v_1)+\cdots+a_mz(v_m).
\end{equation}
Now if $\beta\in [1,N]^k$, then $v(z,\beta)-v\in K_0$,
hence for some $a_{\beta,1},\dots,a_{\beta,m_\beta}\in\Z$,
$v_{\beta,1},\dots,v_{\beta,m_\beta}\in F\cap K_0$, $m_\beta\in\N$:
$$v(\beta,z)-v=a_{\beta,1}v_{\beta,1}+\cdots+a_{\beta,m_\beta}v_{\beta,m_\beta}$$ 
and therefore
\begin{equation}\label{equ2}
z\big(\beta,v(\beta,z)-v\big) = 
a_{\beta,1}z(\beta,v_{\beta,1})+\cdots+a_{\beta,m_\beta}z(\beta,v_{\beta,m_\beta})
\end{equation} with
$z(\beta,v_{\beta,j})\in F_N$ for each $j$. Combining
$$z=z(v) + \sum_{\beta} z\big(\beta,v(z,\beta)-v\big)$$
with \eqref{equ1} and \eqref{equ2} yields an expression of $z$ as a $\Z$-linear
combination of vectors in $F_N$, as required.
\end{proof}

\subsection*{Computing ${\cal H}_{k,\infty}$}
We now specify an algorithm which recursively (in $k$) computes the set
${\cal H}_{k,\infty}$.
The following completion procedure is at the heart of the $k$-th step in the algorithm.
We say that a set $G$ of vector-trees of height $k$  is
{\bf root-complete} if  for all $\cal S,\cal T\in G$ the sum 
$\cal S+\cal T$ has a normal form 
$\cal N$ with respect to $G$ such that $\root(\cal N)=0$.

\begin{algorithm} \label{Completion procedure}
{(Completion procedure)}

\smallskip

\Kw{Input:} \parbox[t]{31em}{a finite 
set $G$ of finite vector-trees
of height $k$.}

\smallskip

\Kw{Output:} \parbox[t]{30em}{a finite set of finite vector-trees of height $k$  which contains $G$ and is root-complete.}

\smallskip

$C:= \{\cal S+\cal T: \cal S,\cal T\in G\}$

\Kw{while} $C\neq \emptyset $ \Kw{do}

\hskip 2.5em ${\cal S}:=$ an element in $C$

\hskip 2.5em $C:=C\setminus\{{\cal S}\}$

\hskip 2.5em ${\cal T}:=\text{normalForm}({\cal S},G)$

\hskip 2.5em \Kw{if} $\root({\cal T})\neq 0$ \Kw{then}

\hskip 5em $G:=G\cup \{{\cal T},-\cal T\}$

\hskip 5em $C:=C\cup\{\cal S+\cal T, \cal S+(-\cal T) : \cal S\in G\}$

\Kw{return} $G$
\end{algorithm}

\noindent
We turn to
termination and correctness of Algorithm~\ref{Completion procedure}:

\begin{prop}\label{Completion procedure is adequate}
Algorithm \ref{Completion procedure} 
terminates and satisfies its specification. 
\end{prop}

\begin{proof}
Let  $G_i$ be the value of $G$ and
$\cal S_i$ be the value of $\cal S$ in the $i$-th pass of the
\Kw{while}-loop in Algorithm~\ref{Completion procedure}, and
$\cal T_i := \text{normalForm}({\cal S}_i,G_i)$.
Suppose that $\root(\cal T_i)\neq 0$ for all $i$ in an infinite subset
$I$ of $\N\setminus\{0\}$.
Now $\cal T\not\red_k\cal T_i$ for all
$\cal T\in G_i$ with $\root(\cal T)\neq 0$, and 
$\cal T_i\in G_j$ for all $i<j$.
In particular $\cal T_i\not\red_k\cal T_j$ for all $i<j$ in $I$.
This contradicts  Corollary~\ref{Sequence of non-reducible trees is finite}.
Hence there is some $n$ such that $\root(\cal T_i)=0$ for all $i\geq n$.
So if $m$ is the size of the set $C$ before the $n$-th iteration of 
the \Kw{while}-loop, then Algorithm \ref{Completion procedure} terminates 
after $m$ more iterations. This shows termination of
Algorithm \ref{Completion procedure}. Correctness is easily shown.
\end{proof}

We say that a set $G$ of vector-trees of height $k$ is
{\bf symmetric} if $\cal S\in G\Rightarrow -\cal S\in G$ for all $\cal S$.
If $G$ is symmetric, then for each $N\geq 1$, 
the subset $\langle G\rangle_N$ of $\Z^{d_{k,N}}$
is also symmetric: $z\in \langle G\rangle_N\Rightarrow -z\in \langle G\rangle_N$ for every $z\in\Z^{d_{k,N}}$.

\begin{remarks}\ 

\begin{enumerate}
\item If the input set $G$ in
Algorithm \ref{Completion procedure} is symmetric, then so is its
output set. If every vector-tree in $G$ has value $0$, then
so does every vector-tree in the output set.
\item For $k=0$ and input set $G=$ a finite symmetric generating set
for the $\Z$-module $\ker(T_0)$, 
Algorithm~\ref{Completion procedure} reduces to
Algorithm~\ref{Completion Algorithm} 
from Section~\ref{Preliminaries: Test Sets} and computes a finite set of
vectors containing $G(T_0)$.
\end{enumerate}
\end{remarks}

Here now is:

\begin{algorithm} \label{Algorithm to compute H_infty}
{(Algorithm to compute ${\cal H}_{k,\infty}$)}

\smallskip

\Kw{Input:} \parbox[t]{31em}{an integer $k\geq 0$.}

\smallskip

\Kw{Output:} \parbox[t]{30em}{a finite symmetric set $G_k$ of finite vector-trees of height $k$ such that $G(A_{k,N})\subseteq \langle G_k\rangle_{N}$ for all $N\geq 1$.}

\smallskip

\Kw{for} $i=0,\dots,k$ \Kw{do}

\hskip 2.5em \Kw{if} $i=0$ \Kw{then}

\hskip 5em   $F_0:=$ a finite symmetric generating set  for $\ker(A_{0,1})$

\hskip 5em  $G:=\big\{\cal T(v) : v\in F_0\big\}$

\hskip 2.5em \Kw{else}

\hskip 5em   $F_i:=$  \parbox[t]{25em}{a finite symmetric generating set  for  $\ker(A_{i,1})$
satisfying the hypothesis of Lemma~\ref{Generators} (for $k=i$)}

\smallskip

\hskip 5em  $\cal T_0 := \cal T\big(\big\{(0,v): v=0\text{ or }
v\in\paths(\cal S) \text{ for some $\cal S\in G_{i-1}$}\big\}\big)$

\hskip 5em  $G:=\{\cal T_0\} \cup \big\{\cal T(v) : v\in F_i\big\}$

\hskip 2.5em $G_i :=$ output of Algorithm~\ref{Completion procedure} applied to $G$
\end{algorithm}

Since
termination of Algorithm~\ref{Algorithm to compute H_infty} follows from
Proposition~\ref{Completion procedure is adequate}, we only need to establish its correctness:

\begin{theorem}\label{Algorithm is ok}
Let $G_0,\dots,G_k$ be the sets computed by
Algorithm~\ref{Algorithm to compute H_infty},
for given input $k$. Then each $G_i$ is a finite symmetric set 
of finite vector-trees of height $i$ 
with $G(A_{i,N})\subseteq\langle G_i\rangle_N$ for all $N\geq 1$.
\textup{(}In particular,
the set consisting of the building blocks of  vector-trees in $G_k$ contains 
${\cal H}_{k,\infty}$.\textup{)}
\end{theorem}

In the proof we use the following immediate consequence of 
Lemma~\ref{Useful Lemma},~(1) and Corollary~\ref{Useful Lemma, Cor}:

\begin{lemma}\label{Normalform 1}
Let $G$ be a root-complete set of vector-trees of height $k$.
Then for every $N\geq 1$ and $y,z\in\langle G\rangle_N$ there exist
$f\in\Z^{d_{k,N}}$  and
$g_1,\dots,g_n\in\langle G\rangle_N$
such that  $f_\epsilon=0$ and
$$y+z=f+\sum_j g_j, \qquad g_j\red y+z\text{ for all $j$.}$$
\end{lemma}

For the proof of Theorem~\ref{Algorithm is ok}, fix integers
$k\geq 0$ and $N\geq 1$. We show, by induction on $i=0,\dots,k$, that
$G(A_{i,N})\subseteq\langle G_i\rangle_N$. 
The case $i=0$ is covered by Remark~(2) following the proof
of Proposition~\ref{Completion procedure is adequate}. Suppose that $i>0$.
By inductive hypothesis, $G_{i-1}$ is
a finite symmetric set 
of finite vector-trees of height $i-1$ 
with $G(A_{i-1,N})\subseteq\langle G_{i-1}\rangle_N$. 
In particular, $\cal T_0$ has value $0$ and
$-\cal T_0=\cal T_0$. Moreover, by the next lemma,
every $f\in\ker(A_{i,N})$ with $f_{\epsilon}=0$ has normal
form $0$ with respect to $\langle \cal T_0\rangle_N$.
In the proof we denote the concatenation
of the finite sequences $\alpha\in [1,N]^s$ and
$\beta\in [1,N]^t$ ($s,t\in\{0,\dots,k\}$) by  $\alpha\beta\in [1,N]^{s+t}$.

\begin{lemma}\label{Normalform 2}
For every $f\in\ker(A_{i,N})$ with $f_{\epsilon}=0$ there are
$h_1,\dots,h_m\in\langle\cal T_0\rangle_N$
such that
$$f=\sum_j  h_j, \qquad \text{$h_j\red f$ for all $j$.}$$
\end{lemma}
\begin{proof}
Let $f_n\in\Z^{d_{i-1,N}}$, $n=1,\dots,N$, be given by
$f_{n,\alpha}:=f_{(n)\alpha}$ for $\alpha\in [1,N]^{s}$, $s=0,\dots,i-1$.
Since  $f\in\ker(A_{i,N})$ and $f_{\epsilon}=0$, we have
$f_1,\dots,f_N\in\ker(A_{i-1,N})$.
Hence for each $n=1,\dots,N$ there are $M(n)\in\N$
and $g_{nm}\in G(A_{i-1,N})$, 
where $m=1,\dots,M(n)$, such that
$$f_n = \sum_m g_{nm} \qquad g_{nm}\red f_n \text{ for all $m$.}$$
Now let $h_{nm}$ be the vector in $\Z^{d_{i,N}}$ defined as follows: for
$\alpha=(\alpha_1,\dots,\alpha_s)\in [1,N]^s$, $s=0,\dots,i$, let
$$h_{nm,\alpha}=\begin{cases}
g_{nm,\alpha'} & \text{if $s>0$ and $\alpha_1=n$} \\
0              & \text{otherwise,}
\end{cases}$$
where $\alpha'=(\alpha_2,\dots,\alpha_s)$. Then $h_{nm}$ can be
constructed from 
$\cal T_0$,
since $G(A_{i-1,N})\subseteq \langle G_{i-1}\rangle_N$ and
$0\in\paths(\cal T_0)$, and
$f=\sum_{n,m} h_{nm}$ with $h_{nm}\red f$ for all $n$ and $m$.
\end{proof}

By Lemma~\ref{Generators} 
and the choice of $F_i$ in
Algorithm~\ref{Algorithm to compute H_infty},
$\langle G_i\rangle_N$ is a symmetric generating set for 
$\ker(A_{i,N})$. Hence
$\langle G_i\rangle_N$ can be used as an input set for the computation of a
Graver basis for the matrix $A_{i,N}$
(Algorithm~\ref{Completion Algorithm}).
By Lemmas~\ref{Normalform 1} and \ref{Normalform 2},
the sum $y+z$ of two elements $y$ and $z$ of $\langle G_i\rangle_N$
has normal form $0$ with respect to $\langle G_i\rangle_N$. Hence
Algorithm~\ref{Completion Algorithm} applied to $A_{i,N}$ and $\langle G_i\rangle_N$
just returns the input set $\langle G_i\rangle_N$; in particular
$G(A_{i,N})\subseteq \langle G_i\rangle_N$ as desired. \qed

\begin{remarks}\ 

\begin{enumerate}
\item Theorem~\ref{Algorithm is ok} gives another proof of Proposition~\ref{Hkinfty is finite}.
\item For $k=1$, the algorithm to compute $\cal H_{k,\infty}$ described above
differs slightly from the one given in \cite{Hemmecke+Schultz:dec2SIP}.
This is because  Algorithm~3.15
in \cite{Hemmecke+Schultz:dec2SIP} is defective: 
to see this, consider (in the notation introduced there)
the pairs $s=g=(0,\{0,1\})$; then $g\red s$ and $s\ominus g=s$, causing 
Algorithm~3.15 to diverge
on the inputs $s$ and $G=\{g\}$.
\end{enumerate}
\end{remarks}

\section{Finding an Optimal Solution}

\noindent
In this final section we outline how the set $G_k$ produced by
Algorithm~\ref{Algorithm to compute H_infty}
can be employed to solve any particular instance 
\begin{equation}\label{Family of multistage problems}
\min\bigl\{c^\intercal z: A_{k,N}z=b,\ z\in\N^{d_{k,N}}\bigr\}
\tag{$\operatorname{IP}_{N,b,c}$}
\end{equation}
of our family of $(k+1)$-stage stochastic integer programs, given
a choice of the number $N$ of scenarios, a right-hand side  
$b\in\Z^{e_{k,N}}$
(where $e_{k,N}:=N^k\cdot l=$ number of rows of the coefficient matrix
$A_{k,N}$), and
 a cost vector $c\in\R^{d_{k,N}}$. 
%
%
%
%
%
%
Throughout we assume that a finite set $G_k$ of
finite vector-trees of height $k$ such that $G(A_{k,N})\subseteq
\langle G_k\rangle_N$ for all $N\geq 1$ (as computed by Algorithm~\ref{Algorithm to compute H_infty})
is at our disposal.

We first concentrate on the problem of augmenting a feasible solution to an
optimal one.
%
%
%
%
%
%
For this we use the recursive algorithm below. We write a vector
$v\in\R^{d_{k,N}}$ (analogous to our practice in the case of integer vectors above) as
$$v=\big(v_\alpha:\alpha\in [1,N]^s, s=0,\dots,k\big),$$ and for
$k>0$ and $i=1,\dots,N$ we put
$$v_{(i)} := \big(v_{(i)\,\alpha}:\alpha\in [1,N]^s, s=0,\dots,k-1\big)
\in\R^{d_{k-1,N}}.$$
Conversely, if $k>0$, given  $r\in\R^{n_k}$ and $N$ vectors
$v_1,\dots,v_N\in\R^{d_{k-1,N}}$,
we denote by $v(r,v_1,\dots,v_N)$
the vector $v\in\R^{d_{k,N}}$ with $v_\epsilon=r$ and
$v_{(i)}=v_i$ for all $i=1,\dots,N$.
With this notation, given $\cal T$ and $N\geq 1$, the set
$\langle\cal T\rangle_N$ consists of all vectors $z\in\Z^{d_{k,N}}$ of the
form $z=v\big(\root(\cal T),v_1,\dots,v_N\big)$, where
$v_i\in\langle\cal T_i\rangle_N$  for
an immediate sub-vector-tree $\cal T_i$ of $\cal T$,
 for each $i=1,\dots,N$.

\begin{algorithm}\label{Most expensive constructible vector}
{(Algorithm to find most expensive constructible vector)}

\smallskip

\Kw{Input:} \parbox[t]{30em}{an integer $N\geq 1$,  a
finite set $G$ of
finite vector-trees of height $k$, a cost vector 
$c\in\R^{d_{k,N}}$, and a vector $z\in\Z^{d_{k,N}}$.}

\smallskip

\Kw{Output:} \parbox[t]{30em}{a vector $v=\operatorname{mostExpensive}(N,G,c,z)$
with $v\in\langle G\rangle_N$ and
$v\leq z$ such that $c^\intercal v$ is maximal with these properties, or
``FAIL'' if no $v\in\langle G\rangle_N$ with
$v\leq z$ exists.}

\smallskip

$G' := \big\{\cal T\in G : \root(\cal T)\leq z_\epsilon\big\}$

\Kw{if} $G'=\emptyset$
\Kw{then} \Kw{return} ``FAIL''

\Kw{while} $G'\neq\emptyset$ \Kw{do}

\hskip 2.5em $\cal T:=$ the element of $G'$ such that $c_\epsilon^\intercal
\root(\cal T)$ is maximal

\hskip 2.5em $G':=G'\setminus\{\cal T\}$

\hskip 2.5em \Kw{if} $k=0$ \Kw{then} \Kw{return} $\root(\cal T)$

\hskip 2.5em $G_{\cal T} := $ the set of immediate sub-vector-trees of $\cal T$

\hskip 2.5em \Kw{for} $i=1$ \Kw{to} $N$ \Kw{do}

\hskip   5em $v_{i} := \operatorname{mostExpensive}(N,G_{\cal T},c_{(i)},z_{(i)})$

\hskip 2.5em \Kw{if} $v_{1},\dots,v_{N}\neq $ ``FAIL'' \Kw{then}
\Kw{return} $v\big(\root(\cal T),v_{1},\dots,v_{N}\big)$

\Kw{return} ``FAIL''

\end{algorithm}

\noindent
Termination and correctness of the procedure above are easily seen.
By the discussion in Section~\ref{Preliminaries: Test Sets} this implies termination and correctness of the following
algorithm:

\begin{algorithm}\label{AugAlg}
(Augmentation algorithm)

\smallskip

\Kw{Input:}  \parbox[t]{30em}{an integer $N\geq 1$, vectors 
$b\in\Z^{e_{k,N}}$, $c\in\R^{d_{k,N}}$,
and a feasible solution $z\in\N^{d_{k,N}}$ to
\eqref{Family of multistage problems}.}

\smallskip

\Kw{Output:} \parbox[t]{30em}{an optimal solution $\operatorname{optimalSolution}(N,b,c,z)$  to
 \eqref{Family of multistage problems}.}

\smallskip

\Kw{while} $v:=\operatorname{mostExpensive}(N,G_k,c,z)
\neq$ ``FAIL'' and $c^\intercal v>0$ \Kw{do}

\hskip 2.5em $z:=z-v$

\Kw{return} $z$
\end{algorithm} 

\noindent
The next algorithm
produces an initial feasible solution from a given solution
(in $\Z^{d_{k,N}}$) to the equation $A_{k,N}z=b$.
Termination and correctness of this procedure follow from
results in \cite{Hemmecke:PSP}; see also Algorithm~\ref{Finding a feasible solution} above.

\begin{algorithm}
(Finding a feasible solution)

\smallskip

\Kw{Input:}  \parbox[t]{30em}{an integer $N\geq 1$, vectors 
$b\in\Z^{e_{k,N}}$, $c\in\R^{d_{k,N}}$,  and a solution $z\in\Z^{d_{k,N}}$ to
\eqref{Family of multistage problems}.}

\smallskip

\Kw{Output:} \parbox[t]{30em}{a feasible solution 
$\operatorname{feasibleSolution}(N,b,c,z)$ to
 \eqref{Family of multistage problems}, or ``FAIL'' if no such
solution exists.}

\smallskip

\Kw{while} $v:=\operatorname{mostExpensive}\big(N,G_k,c(z),z^+\big)
\neq$ ``FAIL'' and $c(z)^\intercal v>0$ \Kw{do}

\hskip 2.5em $z:=z-v$

\Kw{if} $z\geq 0$ \Kw{then} \Kw{return} $z$ \Kw{else} \Kw{return} ``FAIL''

\end{algorithm} 

\noindent
Finally, this leads to our algorithm for solving instances of \eqref{Family of
multistage problems} using $G_k$:

\begin{algorithm}
(Finding an optimal solution)

\smallskip

\Kw{Input:}  \parbox[t]{30em}{an integer $N\geq 1$ and vectors 
$b\in\Z^{e_{k,N}}$, $c\in\R^{d_{k,N}}$.}

\smallskip

\Kw{Output:} \parbox[t]{30em}{an optimal solution to
 \eqref{Family of multistage problems}, or ``FAIL'' if no feasible
solution exists.}

\smallskip

\Kw{if} there is no  $z\in\Z^{d_{N,k}}$ with  $A_{k,N}z= b$ \Kw{then}

\hskip 2.5em \Kw{return} ``FAIL''

\Kw{else}

\hskip 2.5em $z:=$ an element of $\Z^{d_{N,k}}$ satisfying $A_{k,N}z=b$

\smallskip

$f:=\operatorname{feasibleSolution}(N,b,c,z)$

\smallskip

\Kw{if} $f=$ ``FAIL'' \Kw{then}

\hskip 2.5em \Kw{return} ``FAIL''

\Kw{else}

\hskip 2.5em \Kw{return} $\operatorname{optimalSolution}(N,b,c,f)$
\end{algorithm} 

\begin{remark}
The complexity of Algorithm~\ref{Algorithm to compute H_infty} is
unclear. For a very crude complexity result related to
Maclagan's principle see \cite{AschenbrennerPong}, Proposition~3.25.
Some computational experiments in the case $k=1$ are reported in
\cite{Hemmecke+Schultz:dec2SIP}, Section~4.
Note that (once $G_k$ is available) 
the running time of Algorithm~\ref{Most expensive constructible vector} above is $O(N^{k})$.
\end{remark}

\bibliographystyle{amsplain}

\providecommand{\bysame}{\leavevmode\hbox to3em{\hrulefill}\thinspace}
\providecommand{\MR}{\relax\ifhmode\unskip\space\fi MR }
\providecommand{\MRhref}[2]{%
  \href{http://www.ams.org/mathscinet-getitem?mr=#1}{#2}
}
\providecommand{\href}[2]{#2}

\end{document}